\documentclass[11pt
]{article}
\usepackage{amssymb,amsmath,amsthm,amscd}
\usepackage[bbgreekl]{mathbbol}
\usepackage{mathrsfs}
\usepackage{mparhack}

\makeatletter
\if@twoside%
\else\fi
\makeatother

\newcommand{\nek}{\newcommand}
\newcommand{\renek}{\renewcommand}

\nek{\vyk} [1] {}

\vyk{
\addtolength{\topmargin}{-7ex}
\addtolength{\topmargin}{4ex}

\addtolength{\textheight}{-4ex} 
\addtolength{\textwidth}{-3ex}
}

\nek{\noi}{\noindent}

\DeclareMathAlphabet{\cur}{U}{eur}{m}{n}
\DeclareMathAlphabet{\skr}{U}{eus}{m}{n}
\nek{\ubf} {\fontseries{b}\selectfont}

\nek{\bfit}{\bfseries\itshape}%
\nek{\bfsl}{\bfseries\slshape}%

\nek{\parf}{\subsection}
\nek{\punk}{\subsubsection}  
\renewcommand{\thesubsection}{\arabic{subsection}}

\nek{\itla}{\item\label}

\newcounter{enuc}
\nek{\enuci}{\addtocounter{enuc}{1}}
\nek{\fenu}{
\tenu{{\mtho$(\fnsymbol{enumi})$}}\itsep}

\nek{\aenu}{
\tenu{{\rm(\arabic{enumi})}}\itsep}

\nek{\Aenu}{
\tenu{{\rm(\Alph{enumi})}}\itsep}

\nek{\xaenu}{
\tenu{{\rm(\alph{enumi})}}\itsep}

\nek{\renu}{
\tenu{{\rm(\roman{enumi})}}\itsep}

\nek{\Renu}{
\tenu{{\rm(\Roman{enumi})}}\itsep}

\nek{\cenu}{%
\def\theenumi{$\mtho\arabic{enumi}^\circ$}
\itsep}

\nek{\itsep}{\itemsep=0.0ex plus 0.05ex minus 0.1ex}
\nek{\tenu}[1]{
\def\theenumi{#1}
\itsep
}
\nek{\tenui}[1]{
\def\theenumii{#1}
\itsep
}

\nek{\lapm} [1] {{\rm``$#1$''}}
\nek{\thsp}{\hspace{0.1ex}}

\theoremstyle{plain}
\newtheorem{theore}             {Theorem} 
\newtheorem{pro}  [theore]  {Proposition} 
\newtheorem{cor}  [theore]  {Corollary} 
\newtheorem{lem}  [theore]  {Lemma} 
\newtheorem{cla}  [theore]  {Claim} 
\newtheorem{gipot}  [theore]  {Conjecture} 

\theoremstyle{definition}
\newtheorem{opr}  [theore]  {Definition} 
\newtheorem{rem}  [theore]  {Remark} 
\newtheorem{bla}  [theore]  {Blanket agreement} 
\newtheorem{vop}  [theore]  {Problem} 
\newtheorem*{pru}{Proof}
\newtheorem{cru}  [theore]  {Crucial assumption}

\nek{\bvo} {\begin{vop}}
\nek{\evo} {\qed\end{vop}}
\nek{\bcl} {\begin{cla}}
\nek{\ecl} {\qed\end{cla}}
\nek{\bcr} {\begin{cru}}
\nek{\ecr} {\qed\end{cru}}
\nek{\bgi} {\begin{gipot}}
\nek{\egi} {\end{gipot}}
\nek{\bre} {\begin{rem}}
\nek{\ere} {\qed\end{rem}}
\nek{\bte} {\begin{theore}}
\nek{\ete} {\end{theore}}
\nek{\bdf} {\begin{opr}}
\nek{\eDf} {\end{opr}}
\nek{\edf} {\qed\end{opr}}
\nek{\bbl} {\begin{bla}}
\nek{\ebl} {\end{bla}}
\nek{\bco} {\begin{cor}}
\nek{\eco} {\end{cor}}
\nek{\ble} {\begin{lem}}
\nek{\ele} {\end{lem}}
\nek{\bpr} {\begin{pro}}
\nek{\epr} {\end{pro}}

\nek{\bpf} {\begin{pru}}
\nek{\epf} {\qed\end{pru}}
\nek{\epF} [1] {\qed{\,\;({\sl#1\,})}\end{pru}}
\nek{\qeD} [1] {\qed{\;\,({\sl#1\,})}}

\nek{\buss}[1]{\begin{equation}#1\end{equation}}
\nek{\ben}{\begin{enumerate}}
\nek{\een}{\end{enumerate}}
\nek{\bit}{\begin{itemize}}
\nek{\eit}{\end{itemize}}
\nek{\bde}{\begin{description}}
\nek{\ede}{\end{description}}
\nek{\snos} [1] {\,\footnote{\:#1}}%
\nek{\snot} [1] {\footnotetext{\:#1}}%
\nek{\snom}     {\,\footnotemark}%
\nek{\DS}{\displaystyle}
\nek{\etc} {{\sl etc}}
\nek{\noo}
{\text{w.\hspace{0.3ex}l.\hspace{0.3ex}o.\hspace{0.3ex}g.}}
\nek{\ie}  {\hbox{\sl i.\hspace{0.3ex}e.}}
\nek{\vrt}  {\hbox{\sl w.\hspace{0.3ex}r.\hspace{0.3ex}t.}}
\nek{\ea}  {\hbox{\sl et.\hspace{0.3ex}al.}}
\nek{\eg}  {\hbox{\sl e.\hspace{0.3ex}g.}}
\nek{\lsc} {\hbox{l.\hspace{0.3ex}s.\hspace{0.3ex}c.}}
\nek{\qd}  {countably determined}
\nek{\card} {\mathop{\mathtt{card}}}
\nek{\tlim} {\mathop{\mathtt{lim}}}
\nek{\dom}  {\mathop{\text{\tt dom}}}
\nek{\cof}  {\mathop{\mathtt{cof}}}
\nek{\tmin} {\mathop{\mathtt{min}}}
\nek{\tmax} {\mathop{\mathtt{max}}}
\nek{\tsup} {\mathop{\mathtt{sup}}}
\nek{\tinf} {\mathop{\mathtt{inf}}}
\nek{\ran}  {\mathop{\mathtt{ran}}}

\nek{\Ord} {{\mathtt{Ord}}}
\nek{\Ult} {{\mathtt{Ult}}}
\nek{\dd}[1]{$\mtho\hspace{0.2ex}{#1}$-}
\nek{\al} {\alpha}
\nek{\ba} {\beta}
\nek{\ga} {\gamma}
\nek{\Ga} {\varGamma}
\nek{\da} {\delta}
\nek{\ka}{\kappa}
\nek{\la}{\lambda}
\nek{\La}{\varLambda}
\nek{\vt} {\vartheta}
\nek{\ve} {\varepsilon}
\nek{\vpi}{\varphi}
\nek{\Da} {\varDelta}
\nek{\sg} {\sigma}
\nek{\Sg} {\varSigma}
\nek{\za} {\zeta}
\nek{\om} {\omega}
\nek{\omi}{\om_1}

\nek{\dA} {\mathbb A}
\nek{\dD} {\mathbb D}
\nek{\dH} {\mathbb H}
\nek{\dU} {\mathbb U}
\nek{\dP} {\mathbb P}
\nek{\dQ} {\mathbb Q}
\nek{\dT} {\mathbb T}

\nek{\bV} {\mathbf V}
\nek{\bL} {\mathbf L}
\nek{\cF} {\mathscr F}
\nek{\cG} {\mathscr G}
\nek{\cP} {\mathscr P}

\nek{\sd}   {\mathbin{\bigtriangleup}}
\nek{\di}
{\mathop{\text{\mtho\large\boldmath$\Da$}}\limits}

\nek{\ali}{\aleph_1}
\nek{\alo}{{\aleph_0}}

\nek{\dm}  {$$}
\nek{\lra} {\longrightarrow} 
\nek{\imp} {\Longrightarrow} 
\nek{\nin} {\notin}
\nek{\sus} {\exists\hspace{0.3ex}}
\nek{\kaz} {\forall\hspace{0.3ex}}
\nek{\eqv} {\Longleftrightarrow} 
\nek{\ti}  {\times} 
\nek{\sqe} {\sq_{\text{\tt end}}}
\nek{\sue} {\su_{\text{\tt end}}}
\nek{\sq}  {\subseteq}
\nek{\sneq}  {\subsetneqq}
\nek{\su}  {\subset}
\nek{\obr} {^{-1}}
\nek{\iy}  {\infty}
\nek{\mo}  {\models}
\nek{\res} {\restriction}
\nek{\pu}  {\varnothing}
\nek{\bez} {\smallsetminus}
\nek{\fc}
{\mathrel{{{|}{|}\hspace*{-1.1ex}{}-\hspace*{-1.5ex}{}-}}}
\nek{\ans} [1] {\{\hspace{0.01ex}#1\hspace{0.01ex}\}} 
\nek{\ang} [1] {\langle#1\rangle} 
\nek{\stk} [2] {\ang{#1\hspace{0.1ex};\hspace{0.1ex}#2}}
\nek{\sis} [2] {\ans{#1}_{#2}}
\nek{\ens} [2] 
{\{\hspace{0.01ex}#1
\hspace*{0.3ex}{:}\hspace*{0.4ex}
#2\hspace{0.01ex}\}} 

\nek{\va} {\vec a}
\nek{\vc} {\vec c}
\nek{\vu} {\vec u}
\nek{\vv} {\vec v}
\nek{\vx} {\vec x}
\nek{\vy} {\vec y}

\nek{\imar}[1]{\marginpar
[\flushright\footnotesize{#1}]%
{\flushleft\footnotesize{#1}}%
}

\nek{\pfi} {^{\text{\rm fin}}}

\nek{\mem} {\dd\in}

\nek{\spp} [3] {{\text{\tt bas}}_{#2}(#3)}
\nek{\jt } [3] {{J_{#2}(#3)}}
\nek{\spi} [1] {{\text{\tt bas}}_{#1}}
\nek{\ppi} [1] {{\text{\tt bas}}'_{#1}}

\nek{\zo} {,\linebreak[0]}
\nek{\zi} {,\linebreak[0]\hspace*{0.2ex}}
\nek{\zd} {,\linebreak[0]\:}
\nek{\zt} {,\linebreak[0]\;}
\nek{\zz} {\linebreak[0]}

\nek{\ema}{\ensuremath}
\nek{\mast} {{\mtho\ema\ast}}
\nek{\trm} {\underline}
\nek{\td} {{\trm d}}
\nek{\uh} {{\trm \zC}}
\nek{\uK} {{\trm K}}
\nek{\tG} {{\trm G}}
\nek{\tX} {{\trm X}}

\nek{\wh} {\widetilde h}
\nek{\ta}  {\check a}
\nek{\tu}  {\check u}
\nek{\tv}  {\check v}
\nek{\tw}  {\check w}
\nek{\tx}  {\check x}
\nek{\ty}  {\check y}
\nek{\tal} {\check\al}
\nek{\tba} {\check\ba}
\nek{\tom} {\check\om}
\nek{\tga} {\check\ga}
\nek{\tda} {\check\da}
\nek{\tka} {\check\ka}
\nek{\tnu} {\check\nu}
\nek{\tet} {\check\eta}
\nek{\tvt} {\check\vt}
\nek{\txi} {\check\xi}
\nek{\tdP} {\check\dP}    
\nek{\tE}  {\check E}    
\nek{\tF}  {\check F}    
\nek{\tI}  {\check I}    
\nek{\tK}  {\check K}    
\nek{\tP}  {\check P}    
\nek{\tQ}  {\check Q}    
\nek{\tR}  {\check R}    
\nek{\tS}  {\check S}    
\nek{\tU}  {\check U}    
\nek{\tY}  {\check Y}    

\nek{\ev} {\equiv_\bV} 
\nek{\lv} {\le_\bV}

\renek{\le} {\leqslant}

\nek{\jj} [2] {
{{#1}\mathbin{{{/}\hspace*{-0.7ex}{/}\hspace*{-0.2ex}}}{#2}}}

\nek{\nx} [1] {{#1}^{\dagger}}

\nek{\gc}{{\mathfrak c}}
\nek{\cX}{{\mathscr X}}

\nek{\mtho}{\mathsurround=0mm}
\nek{\msur}{\hspace*{-1\mathsurround}}
\nek{\vom} {\vspace{0.5ex}}
\nek{\vtm} {\vspace{1ex}}

\nek{\bM} {\text{\bfit M}}
\nek{\hM} {\bM_\om}
\nek{\ha} {\ka_\om}
\nek{\bN} {\text{\bfit N}}
\nek{\les} {\leqslant}
\nek{\ges} {\geqslant}

\nek{\rit} [1] {{\it#1\/}}

\nek{\zC} {h}

\nek{\spo} [1] {\spi{#1}}

\nek{\jjo} [2] {\jj{#1}{\spi{#2}}}

\nek{\fo} [1] {F({#1})}

\nek{\fa} [1] {F(\ans{#1})}

\nek{\nn} {n}

\nek{\wa} {t} 

\nek{\evn} [1] {\jj{#1}{\EV}}
\nek{\EV} {{\text{\sc ess}}}
\nek{\eva} {\EV^\ast}

\nek{\uhev} {\underline d}
\nek{\hev} {\evn\zC}

\nek{\xip}{{\xi+1}}
\nek{\kap}{\ka^+}
\nek{\kapn}[1]{\ka^{+#1}}
\nek{\kao} {\kapn\om}

\nek{\sld} [1] {{#1}^\oplus}
\nek{\prd} [1] {{#1}^\ominus}

\nek{\cK} {\mathscr K}

\nek{\onto} {\overset{\text{\tt onto}}\longrightarrow}

\nek{\dpd} {\dP/d}

\renek{\eva} {\EV}

\nek{\ale} {\aleph}
\nek{\alom} {\ale_\om}

\nek{\aal} [1] {[\ale_{#1},\ale_{#1+1})}
\nek{\oal} [1] {[\om,\ale_{#1+1})}
\nek{\aad} [1] {\aal{#1}^2}
\nek{\oad} [1] {\oal{#1}^2}

\nek{\omal} {[\om,\alom)}

\nek{\xd} {D^{(2)}}
\nek{\bxd} {\bd^{(2)}}

\nek{\xP} {\dP^{(2)}}
\nek{\xpf} {\xP_{\text{\tt fin}}}

\nek{\tq} {q_{{\ast}}}
\nek{\tp} {p_{{\ast}}}
\nek{\tr} {r_{{\ast}}}
\nek{\tiv} {v_{{\ast}}}
\nek{\tiu} {u_{{\ast}}}
\nek{\tpa} {{p'}_{{\ast}}}

\nek{\eq} {e^q}
\nek{\ep} {e^p}
\nek{\epa} {e^{p'}}

\nek{\abs} [1] {|#1|}

\nek{\seq} [2] {(#1)_{#2}}
\nek{\sen} [3] {\seq{{#1}\og{#2}}{#3}}

\nek{\pil} [1] {\Pi^{\le#1}}

\nek{\bd} {\overline D}

\nek{\ap}{\cdot}

\nek{\bay}{\begin{array}}
\nek{\eay}{\end{array}}

\nek{\lam} [1]
{\label{#1}\hspace*{0pt}\imar{#1}}%
\nek{\las} [1]
{\label{#1}\imae{#1}}%

\nek{\imae}[1]{\marginpar
[\flushright\footnotesize{#1}]%
{\flushleft\footnotesize{#1}}
}%

\nek{\ogr} [1] {\hspace*{0.2ex}[#1]}
\nek{\og} {\ogr}
\nek{\ogl} [1] {\hspace*{0.2ex}[\le#1]}

\nek{\ogd} [1] {\dD\og{#1}}
\nek{\ogh} [1] {\dH\og{#1}}
\nek{\oph} [1] {\pH\og{#1}}
\nek{\ogp} [1] {\dP\og{#1}}
\nek{\opp} [1] {\pP\og{#1}}

\nek{\pP} {\dP^+}
\nek{\pH} {\dH^+}
\nek{\pQ} {\dQ^+}
\nek{\pqa} {\dQ^\ast}

\nek{\dtr} {\dT^{\text{\tt uni}}}

\nek{\Psd} {\mathbb\Psi}
\nek{\Phd} {\mathbb\Phi}

\nek{\odp} [1] {\xP\og{#1}}

\nek{\re} [2] {{#1}{\res}_{\le#2}}
\nek{\reo} [2] {{\ogp{#1}}{\res}_{\le#2}}
\nek{\rep} [1] {{\dP}{\res}_{\le#1}}
\nek{\req} [1] {{\dQ}_{#1}}
\nek{\qle} [1] {{\dQ}{\res}_{\le#1}}
\nek{\tle} [1] {\dT_{\le#1}}

\nek{\df} {D_{\text{\tt fin}}}
\nek{\xdf} {D^{(2)}_{\text{\tt fin}}}

\nek{\bas}  {\mathop{\text{\tt bas}}}
\nek{\Dom}  {\mathop{\texttt{DOM}}}

\nek{\oa} {\bar a}
\nek{\ove} {\bar e}

\nek{\vF} {\vec F}
\nek{\vX} {\vec X}
\nek{\vY} {\vec Y}

\nek{\vza} {\vec\za}
\nek{\vvt} {\vec\vt}

\nek{\Vs} {\bV_{\text{\tt sym}}}
\nek{\Ls} {\bL_{\text{\tt sym}}}

\nek{\ima} [2] {#1[#2]}

\nek{\dqs} {\dQ^{\text{\tt sat}}}

\nek{\pif} {\Pi_{\text{\tt fin}}}

\nek{\vr} {{\vec \rho}}
\nek{\vta} {{\vec \tau}}

\nek{\app}{\circ}

\nek{\bap} {\bar p}
\nek{\baq} {\bar q}

\nek{\tip} {\tilde p}
\nek{\tiq} {\tilde q}

\nek{\ddom} [1] {\dom{(\dom{(#1)})}}

\nek{\rbn} [2] {\text{\tt Ran}_{#1}(#2)}

\nek{\pes} 
{{\hspace*{0.1ex}{\restriction}\hspace*{-0.65ex}{\restriction}\hspace*{0.1ex}}}

\nek{\mf} [1] {\Phi_{#1}}

\nek{\nd} [1] {\dD\og{#1}}

\nek{\itlm} [2] {\itla{#1}#2\imar{#1}}

\nek{\atc} [1] {\addtocounter{enumi}{#1}}

\renek{\ens} [2] 
{\{\hspace{0.01ex}#1
\hspace*{0.4ex}{:}\hspace*{0.5ex}
#2\hspace{0.01ex}\}} 

\nek{\ftp} [2] {{#1}\downarrow{#2}}
\nek{\ftw} [2] {{#1}\vec{\phantom{x}}{#2}}

\nek{\equ} {=^*}

\nek{\dda} {\dD^\ast}

\renek{\res} 
{{\hspace*{0.15ex}{\restriction}\hspace*{0.15ex}}}

\nek{\mto} {\longmapsto}

\nek{\txy} {\dT_{N\Ga}[G]}

\renek{\punk} [1] {\subsubsection{\ubf #1}}

\nek{\zza} {{\mathbf a}}  
\nek{\zzf} {{\mathbf x}}  
\nek{\zzy} {{\mathbf y}}  
\nek{\zzz} {{\mathbf r}}  
\nek{\zzF} {{\mathbf X}}  
\nek{\zzY} {{\mathbf Y}}  
\nek{\vzF} {\vec\zzF}  
\nek{\vzY} {\vec\zzY}  

\nek{\Fg} {\zzF^G}
\nek{\Fh} {\zzF^H}
\nek{\fg} {\zzf^G}
\nek{\fh} {\zzf^H}
\nek{\yg} {\zzy^G}
\nek{\yh} {\zzy^H}
\nek{\Yg} {\zzY^G}
\nek{\Yh} {\zzY^H}
\nek{\ag} {\zza^G}
\nek{\ah} {\zza^H}
\nek{\wg} {W[G]}
\nek{\wH} {W[H]}
\nek{\vyg} {\vzY[G]}
\nek{\vfg} {\vzF[G]}
\nek{\vyh} {\vzY[H]}
\nek{\vfh} {\vzF[H]}

\nek{\vf} {\vec \zzf}
\nek{\vz} {\vec \zzz}
\renek{\vy} {\vec \zzy}

\nek{\sw} [2] {\mathbf S_{#1#2}}
\nek{\swt} [3] {{\mathbf s}_{#1#2}^{#3}}

\nek{\qs} {\supseteq}

\nek{\psur}{\hspace*{1\mathsurround}}

\nek{\eqn} [1] {\equiv_{#1}}
\nek{\eqa} {\equiv^{\ast}}

\nek{\vxng} {\vec{\mathbf x}_N[G]}
\nek{\vyng} {\vec{\mathbf y}_\Ga[G]}
\nek{\vxnh} {\vec{\mathbf x}_N[H]}
\nek{\vynh} {\vec{\mathbf y}_\Ga[H]}

\renek{\parf}{\section}
\renek{\punk}{\subsection}  
\renek{\thesection}{\arabic{section}}
\renek{\thesubsection}{\arabic{subsection}}

\renek{\lam} [1] {\label{#1}}
\renek{\las} [1] {\label{#1}}
\renek{\imar} [1] {}

\begin{document}

\title{On automorphisms behind the Gitik -- Koepke model 
for violation of the Singular Cardinals Hypothesis w/o 
large cardinals}

\author{Vladimir Kanovei}

\date{\today}

\maketitle

\begin{abstract}
It is known that the assumption that 
``GCH first fails at $\aleph_\om$'' 
leads to large cardinals in {\bf ZFC}. 
Gitik and Koepke~\cite{gk}
demonstrated that this is not so in {\bf ZF}: 
namely there is a generic cardinal-preserving extension of 
$\bL$ (or any universe of {\bf ZFC} + GCH) in which all 
{\bf ZF} axioms hold, the axiom of choice fails, 
$\card{2^{\aleph_n}}=\aleph_{n+1}$ for all natural $n$, 
but there is a surjection from $2^{\alom}$ onto $\la$, 
where $\la>\aleph_{\om+1}$ is any previously chosen 
cardinal in $\bL$, for instance, $\aleph_{\om+17}$.
In other words, in such an extension GCH holds in proper 
sense for all cardinals $\aleph_n$ but fails at $\alom$ 
in Hartogs' sense. 

The goal of this note is to analyse the system of automorphisms 
involved in the Gitik -- Koepke construction.
\end{abstract}

It is known (see \cite{cite1}) that the consistency of the statement 
``GCH first fails at $\alom$'' with {\bf ZFC} definitely requires 
a large cardinal. 
Gitik and Koepke~\cite{gk}
demonstrated that picture changes in the absense of the axiom 
of choice, if one agrees to treat the 
violation of GCH in Hartogs' sense.  
Namely there is a generic cardinal-preserving extension of 
$\bL$ (or any universe of {\bf ZFC} + GCH) in which all 
{\bf ZF} axioms hold, the axiom of choice fails, 
$\card{2^{\aleph_n}}=\aleph_{n+1}$ for all natural $n$, 
but there is a surjection from $2^{\alom}$ onto $\la$, 
where $\la>\aleph_{\om+1}$ is any previously chosen 
cardinal in $\bL$, for instance, $\aleph_{\om+17}$.
Thus in such an extension GCH holds in proper 
sense for all cardinals $\aleph_n$ but fails at $\alom$ 
in Hartogs' sense. 

For the sake of convenience we formulate the main result 
as follows.

\bte
[Gitik -- Koepke \cite{gk}]
\lam M
Let\/ $\la>\aleph_{\om+1}$ be a cardinal in\/ $\bL$, 
the constructible universe.
There is a set-generic extension\/ $\bL[G]$ of\/ $\bL$  
and a symmetric cardinal-preserving 
subextension\/ $\Ls[G]\sq\bL[G]$, such that 
the following is true in\/ $\Ls[G]$$:$
\ben
\renu
\itla{M1}
all axioms of\/ {\bf ZF}$;$

\itla{M2}\msur
$\card{2^{\aleph_n}}=\aleph_{n+1}$ for all natural\/ $n$$;$

\itla{M3}
there is a surjection from\/ $2^{\alom}$ onto\/ $\la$.
\een 
\ete

The goal of this note is to analyse the system of automorphisms 
(which turns out to consist of three different subsistems)
involved in the Gitik -- Koepke proof of this theorem 
in \cite{gk}.\snos
{The author learned the description of the Gitik -- Koepke 
model in the course of his visit to Bonn in the Winter of 
2009/2010.}
On the base of our analysis, we present the proof in a 
somewhat more pedestrian way than in \cite{gk}.

\parf{Basic definitions and the forcing}
\las{bd}
\las{bdf}

After an array of auxiliary definitions, we'll introduce 
the forcing.
 
$\la$ is a fixed cardinal everywhere; $\la>\aleph_\om$.

\punk{Basic definitions} 
\las{bd1}

We define:
\bde
\item[$\ogd n$ =]
all sets $d\sq\aal n$ such that $\card d\le\aleph_n$
\imar{ogd n}

\item[$\opp n$ =] all functions $p:\dom p\to 2$,
such that $\pu\ne\dom p\sq\aal n$,
\imar{opp n}

\item[$\ogp n$ =] all functions $p\in\opp n$,
such that $\dom p\in \ogd n$,
\imar{ogp n}

\item[$\dD$ =]
all sets $d\sq\omal$ such that $d\cap\aal n\in\dD\og n$
for all $n$,

\item[$\dda$ =]
all sets $d\sq\omal$ such that there is $n_0\in\om$ 
such that $d\cap\aal n\in\dD\og n$ for all $n\ge n_0$,

\item[$\pP$ =] all functions $p:\dom p\to 2$ such that 
$\dom p\sq\omal$,
\imar{pP}

\item[$\dP$ =] all functions $p\in\pP$ such that 
$\dom p\in\dD$.
\imar{dP}
\ede
If $n\in\om$ then we let $d\og n=d\cap\aal n$ and 
$p\og n=p\res{\aal n}$ for all $d\in\dD$ and $p\in\pP$. 
Thus $d\in\dD$ iff $d\og n\in\ogd n$ for all $n$, and 
$p\in\dP$ iff $p\og n\in\ogp n$ for all $n$.

We order $\dP$ so that $p\le q$ iff $\dom q\sq\dom p$ and 
$q=p\res{\dom q}$.

Note that if $m\ne n$ then $\ogp n\cap\ogp m=\pu$.

\punk{Assignments} 
\las{bd2}

An \rit{assignment} will be any function $a$ such that 
\ben
\tenu{{\rm(a\arabic{enumi})}}
\itla{a1}\msur 
$\dom a=\bas a\ti\abs a$, where 
$\bas a\sq\om$ and $\abs a\sq\la$ are finite sets, and 

\itla{a2}
if $\ang{n,\ga}\in\dom a$ then $a(n,\ga)\in\aal n$.
\een
In particular, $\pu$ (the empty assignment) belongs to 
$\dA$.\snos 
{We suppose that $\bas\pu=\abs\pu=\pu$, but it can be 
consistently assumed that either $\bas\pu=\pu$ and 
$\abs \pu=\Ga\sq\la$ is any finite set, or $\abs\pu=\pu$ and 
$\bas\pu=N\sq\om$ is any finite set, depending on the context. 
Any assignment $a\ne\pu$ has definite 
values of $\abs a$ and $\bas a$.}

If $n\in\bas a$ then define a map $a\og n$   
on the set $\abs a$ by $a\og n(\ga)=a(n,\ga)$.

The set $\dA$ of all assignments is \rit{ordered} so that 
$a\le b$ ($a$ is stronger) iff
\ben
\tenu{{\rm(a\arabic{enumi})}}
\addtocounter{enumi}2
\itlm{a3}\msur 
$\bas b\sq\bas a$ and $\abs b\sq\abs a$, and 

\itlm{a4}
if $n\in\bas a\bez\bas b$ and $\ga\ne\da$ 
belong to $\abs b$ then $a(n,\ga)\ne a(n,\da)$.
\een
Clearly $\pu$ is the \dd\le largest element in $\dA$.

Assignments $a,b$ are \rit{coherent} iff $\dom a=\dom b$, 
and for any $n\in\bas a=\bas b$ and $\ga,\da\in\abs a=\abs b$ 
 we have: 
$a(n,\ga)=a(n,\da)$ iff $b(n,\ga)=b(n,\da)$. 

If $a\in\dA$ and $\Da\sq \abs a$ then let
$a\pes\Da$ be the restriction $a\res{(\bas a\ti\Da)}$.


\punk{Narrow subconditions} 
\las{bd3}

Let $\pH$ consist of all indexed sets 
\imar{pH}
$h=\sis{h_\xi}{\xi\in\abs h}$, where
$\abs h\sq\omal$ and 
$h_\xi\in\opp n$ for all $n$ and $\xi\in\abs h\cap\aal n$. 

We put $h\og n=h\res{\aal n}$ (restriction) 
for $h\in\pH$ and any $n$. 
Thus still $h\og n\in\pH$ and $\abs{h\og n}=\abs h\cap\aal n$. 

Let $\dH$ consist of all   
\imar{dH}
$h\in\pH$ such that 
\ben
\tenu{{\rm(h\arabic{enumi})}}
\itlm{h1}\msur
$\card\abs{h\og n}\le \aal n$ for all $n$,

\itlm{h2}
the set $\bas h=\ens{n}{h\og n\ne \pu}$ is finite,

\itlm{h3}\msur
$h_\xi\in\ogp n$ for all $n$ and $\xi\in\abs h\cap\aal n$. 
\een

We say that a condition $h\in\dH$ is   
\bde
\itsep
\item[\rit{regular}] 
at some $n\in\bas h$, iff for every 
$\xi\in\abs{h}\cap\aal n$ the set 
$\ens{\eta\in\abs{h}\cap\aal n}{h_\eta=h_\xi}$ 
has cardinality exactly $\aleph_n$, 

\item[\rit{stronger}] 
than another condition $g\in\dH$,  
symbolically $h\le g$, iff $\abs g\sq\abs h$, 
and $h_\xi\le g_\xi$ for all $\xi\in\abs g$.
\ede
The empty condition $\pu\in\dH$ ($\abs \pu=\pu$) is 
\dd\le largest in $\dH$.

We further define
$\ogh n=\ens{h\in\dH}{\abs h\sq\aal n}$; 
\imar{ogh n}%
thus $\ogh n$ consists of all indexed sets 
$h=\sis{h_\xi}{\xi\in\abs h}$, where $\abs h\in\dD\og n$
(that is, $\abs h\sq\aal n$ and $\card{\abs h}\le\aleph_n$), 
and $h_\xi\in\ogp n$ for all $\xi\in\abs h$.

It is clear that $h\in\dH$ iff $h\og n\in\ogh n$ for all 
$n$ and the set $\bas h$ is finite.

\punk{Wide subconditions} 
\las{bd4}

Let $\pQ$ consist of all indexed sets 
$q=\sis{q_\ga}{\ga\in\abs q}$, where
$\abs q\sq\la$ and $q_\ga\in\pP$ 
for all $\ga\in\abs q$.
We define
\bde
\item[$\pqa$ =]
all $q\in\pQ$ such that
$\abs q$ is finite,
\imar{pqa}

\item[$\dQ$ =] 
all $q\in\pQ$ such that
$\abs q$ is finite and $q_\ga\in\dP$ 
for all $\ga\in\abs q$.
\imar{dQ}
\ede
 
We say that a condition $q\in\pQ$ is: 
\bde
\itsep
\item[\rit{uniform}\rm,] 
if 
$\dom q_\ga\og n=\dom q_\da\og n$ 
for all $\ga,\da\in\abs q$ and $n\in\om$, 

\item[\rit{compatible}] 
with an assignment $a\in \dA$, iff 
we have $q_\ga\og n=q_\da\og n$ 
whenever $\ga,\da\in\abs q\cap\abs a$, $n\in\bas a$, and 
$a(n,\ga)=a(n,\da)$.

\item[\rit{equally shaped}] 
with another condition $p\in \pQ$, iff
$\abs p=\abs q$, and we have $\dom{p_\ga\og n}=\dom{q_\ga\og n}$ 
holds for all $\ga\in \abs p$ 
and $n\in\om$.

\item[\rit{stronger}] 
than another condition $p\in \pQ$,  
symbolically $q\le p$, iff $\abs p\sq\abs q$, 
and $p_\ga\le q_\ga$ in $\dP$ for all $\ga\in\abs p$.
\ede
Once again, the empty condition $\pu\in\dQ$ ($\abs \pu=\pu$) is 
\dd\le largest in $\dQ$.

\punk{Conditions} 
\las{bd5}

Let $\dT$, {\ubf the forcing}, consist of all triples 
of the form 
$t=\ang{q^t,a^t,h^t}$, where 
$q^t\in\dQ$, $a^t\in\dA$, $h^t\in\dH$, and
\ben
\tenu{{\rm(t\arabic{enumi})}}
\itlm{t1}\msur
$\abs{a^t}=\abs{q^t}$ and $\bas{a^t}=\bas{h^t}$ 
--- we put $\abs{t}:=\abs{a^t}$ 
and $\bas{t}:=\bas{a^t}$,

\itlm{t3}\msur
$\ran a^t\sq\abs{h^t}$ and we have $h^t_{a^t(n,\ga)}=q^t_\ga\ogr n$ 
for all $n\in\bas t$ and $\ga\in \abs t$.

\itlm{t2} 
therefore
$q^t$ is compatible with $a^t$ in the sense above, that is, 
if $\ga,\da\in\abs t$, $n\in\bas t$, and 
$a^t(n,\ga)=a^t(n,\da)$, then 
$q^t_\ga\og n=q^t_\da\og n$.
\een
The set $\dT$ is ordered componentwise: 
a condition $t\in\dT$ is \rit{stronger than} $s\in\dT$, 
symbolically $t\le s$, iff 
$q^t\le q^s$ in $\dQ$, $a^t\le a^s$ in $\dA$, $h^t\le h^s$ in $\dH$. 
Clearly $t=\ang{\pu,\pu,\pu}$ is the largest condition in $\dT$.

A condition $t\in\dT$ is \rit{uniform}, symbolically 
$t\in\dtr$, iff $q^t$ is uniform.

\vyk{
\bdf
\lam{esha}
Conditions $p,q\in \dQ$ are \rit{equally shaped} iff
$\bas p=\bas q$ and $\abs p=\abs q$, 
and \rit{strongly equally shaped} iff in addition $a^p=a^q$.
\edf 

\bdf
[restrictions]
\lam{rest}
Let $\Da\sq\la$ be a finite set. 

If $a\in\dA$ and $\Da\sq \abs a$ then let
$a\pes\Da$ be the restriction $a\res{(\bas a\ti\Da)}$. 

If $q\in\dQ$ and $\Da\sq\abs q$ then let 
$q\pes\Da=\sis{q_\ga}{\ga\in\Da}$.\snos 
{This is not a condition in $\dQ$ since it has no $a^q$.}  
\edf

{\ubf Comment: the forcing as a product.}  
For any $n$ let $\dQ\og n$ consist of all indexed sets 
$q=\ang{a^q,\sis{q_\ga}{\ga\in\abs q}}$, where
\ben
\tenu{{\rm(p\arabic{enumi})}}
\itlm{31}\msur
$\abs q\sq\la$ is finite and $q_\ga\in\ogp n$ 
for all $\ga\in\abs q$,

\itlm{33}\msur
$a^q$ is a function from $\abs q$ to $\aal n$,

\itlm{36}
if $\ga,\da\in\abs q$ and $a^q(\ga)=a^q(\da)$  
then $q_\ga=q_\da$.
\een
Then $\dQ$ as a whole can be identified with the product 
$\prod_{n=0}^\om\dQ\og n$ defined so that it is a finite 
support product \vrt\ the components $a^q$, 
and unrestricted (countable) 
support product \vrt\ the components $q_\ga$. 
}

\parf{Permutations}
\las{sym1}

In this section and the following two sections we 
consider three groups of full or partial order-preserving 
transformations of conditions. 


Let $\pif$ be the group of all permutations of 
the set $\omal$
such that
\ben
\tenu{(\Alph{enumi})}
\itlm{p1} 
for any $n$, the restriction $\pi\og n=\pi\res\aal n$ 
is a permutation of the set $\aal n$,

\itlm{p2} 
the set $\bas\pi=\ens{n}{\pi\og n\ne\text{ the identity}}$ 
is finite.
\een
Let $\pif\og n$ consist of all $\pi\in\pif$ equal to the identity 
outside of $\aal n$. 
Any $\pi\in\pif\og n$ is naturally identified with $\pi\og n$.

There are two types of induced action of 
transformations $\pi\in \pif$, namely:
\ben
\Renu
\itlm{api1}
if $f$ is a function such that $\ran f\sq\omal$ then 
$f'=\pi\ap f$ is a function with the same domain and 
$f'(x)=\pi(f(x))$ for all $x\in\dom f=\dom f'$;

\itlm{api2}
if $f$ is a function such that $\dom f\sq\omal$ then 
$f'=\pi\ap f$ is a function,   
$\dom f'=\ens{\pi(\xi)}{\xi\in\dom f}$, and $f'(\pi(x))=f(x)$ 
for all $\xi\in\dom f$.\snos
{We ignore the conflicting case when both $\ran f\sq\omal$ and 
$\dom f\sq\omal$ as it will never happen in the domains of 
action of transformations $\pi\in \pif$ considered below.} 
\een

Accordingly, we define that any $\pi\in \pif$:
\ben
\aenu
\itlm{pe1} 
acts on $\dA$ by \ref{api1}, so that if 
$a\in\dA$ then $a'=\pi\ap a\in\dA$, $\dom a'=\dom a$, and 
$a'(n,\ga)=\pi(a(n,\ga))$ for all $\ang{n,\ga}\in\dom a$;

\itlm{pe2} 
acts on $\pH$ (and on $\dH\sq\pH$) by \ref{api2}, so that if 
$h\in\pH$ then $h'=\pi\ap h\in\pH$, 
$\abs{h'}=\ens{\pi(\xi)}{\xi\in\abs h}$, and 
$h'_{\pi(\xi)}=h_\xi$ for all $\xi\in\abs h$.
\vyk{
\itlm{pe3}
acts on $\dQ$ as the identity, $\pi\ap q=q$ for any $q\in\dQ$, 
since neither \ref{api1} nor 
\ref{api2} is, generally speaking, applicable to an arbitrary 
element $q\in\dQ$.
}%
\een

Finally if $t=\ang{q^t,a^t,h^t}\in \dT$ then put
$\pi\ap t=\ang{q^t,\pi\ap a^t,\pi\ap h^t}$.

The following lemma is rather obvious.

\ble
\lam{sk1}
Any\/ $\pi\in\pif$ is an order-preserving automorphism of the 
ordered sets\/ $\dA$, $\dH$, and\/ $\dT$. 
Moreover if\/ $a\in\dA$ and\/ $n\in\bas a\bez\bas \pi$ then\/ 
$(\pi\ap a)\og n=a\og n$, and accordingly if\/ $h\in\dH$ and\/ 
$n\nin\bas \pi$ then\/ $(\pi\ap h)\og n=h\og n$.\qed
\ele

\parf{Swaps}
\las{sym2}

Suppose that $a,b\in\dA$, 
$\dom a=\dom b=D$, and $\ran a=\ran b$.
Such a pair of assignments induces a 
\rit{swap transformation} $\sw ab$,  
acting: 
$$
\bay{clccc}
\text{from}&
\dA_{a}=\ens{c\in\dA}{c\le a}&
\text{to}&\dA_b\,,& \\[1ex]
\text{from}&
\pQ_{a}=\ens{q\in\pQ}
{\abs a\sq\abs q\land q\text{ is compatible with } a}&
\text{to}&\pQ_b\,,\\[1ex]
\text{from}&
\dQ_{a}=\ens{q\in\dQ}
{\abs a\sq\abs q\land q\text{ is compatible with } a}&
\text{to}&\dQ_b\,.
\eay
$$
Recall that $q\in\pQ$ is compatible with $a\in \dA$ iff 
$q_\ga\og n=q_\da\og n$ holds
whenever $\ga,\da\in\abs a\cap\abs q$, $n\in\bas a$, and 
$a(n,\ga)=a(n,\da)$. 
Obviously $\dQ_{a}=\pQ_a\cap\dQ$.

The action of $\sw ab$ on $\dA_{a}$ is defined as follows: 
\ben
\aenu
\itlm{sabA}
if $c\in\dA_a$ then $c'=\sw ab\ap c\in \dA$, $\dom c'=\dom c$, 
${c'\res D}=b$ (where $D=\dom a=\dom b$), and 
$c'\res{(\dom c\bez D)}= c\res{(\dom c\bez D)}$.
\een

The action of $\sw ab$ on $\pQ_{a}$ is defined as follows. 
First of all, if $n\in\bas a$ and $\ga\in\abs a$ then let 
$\swt abn(\ga)$ be the least $\vt\in\abs a$ satisfying 
$a(n,\vt)=b(n,\ga)$; such ordinals $\vt$ exist because 
$\ran a=\ran b$. 
Thus $\swt abn:\abs a\to\abs a$.
Then: 
\ben
\aenu
\atc1
\itlm{sQ} 
if $q\in\pQ_a$ then $q'=\sw ab\ap q\in \pQ$, $\abs{q'}=\abs q$, 
and for all $n\in\om$ and $\ga\in\abs q$:

\ben
\itlm{sQ1} 
if $\ga\in\abs{a}$ and $n\in\bas{a}$ 
then $q'_\ga\og n=q_\vt\og n$, where $\vt=\swt abn(\ga)$, 

\itlm{sQ2} 
if either $\ga\nin\abs{a}$ or $n\nin\bas{a}$ 
then $q'_\ga\og n=q_\ga\og n$.
\een
\een

Finally if $t\in \dT_a=\ens{t\in\dT}{a^t\le a}$ 
(then $a^t\in\dA_a$ and $q^t\in\dQ_a$) 
then put
$$
\sw{a}{b}\ap t=\ang{\sw{a}{b}\ap q^t,\sw{a}{b}\ap a^t,h^t}.
$$

\ble
\lam{sk2}
Assume that\/ $a,b\in\dA$, $\bas a=\bas b=B$, 
$\abs a=\abs b=\Da$, and\/ $\ran a=\ran b$. 
Then\/ $\sw ab$ is an order-preserving 
bijection\/ 
$\dA_{a}\onto\dA_{b}$, $\dQ_{a}\onto\dQ_{b}$,  
$\dT_{a}\onto\dT_{b}$ and\/ 
$\sw ba$ is the inverse in each of the three cases. 

Lett\/ $t\in\dT_a$. 
Then\/ $t'=\sw ab\ap t\in\dT_b$, $\abs t=\abs{t'}$, 
$\bas t=\bas{t'}$, and$:$
\ben
\renu
\itlm{sk2i} 
if\/ $t$ is uniform, then so is\/ $t'$ and\/ 
$q^t,q^{t'}$ are equally shaped$;$

\itlm{sk2ii} 
if\/ $n\in B$, $\ga\in\Da$, and\/ $a(n,\ga)=b(n,\ga)$ then\/ 
$a^t(n,\ga)=a^{t'}(n,\ga)=a(n,\ga)=b(n,\ga)$ and\/ 
$q^{t'}_\ga\og n=q^t_\ga\og n\,;$

\itlm{sk2iii} 
if\/ $n\in\abs t$ then\/ 
$\ens{q^{t'}_\ga\og n}{\ga\in\abs{t'}}=
\ens{q^{t}_\ga\og n}{\ga\in\abs{t}}\,.$
\een
\ele 
\bpf
The first essential part of the lemma is to show that 
if $t\in \dT_a$ then $t'=\sw ab\ap t\in\dT_b$. 
Basically it's enough to show that $t'\in\dT$. 
And here the only notable task is to prove \ref{t3} of 
Section \ref{bd5}, that is, 
$q^{t'}_\ga\ogr n=h^{t'}_{a^{t'}(n,\ga)}$ for all $n\in\bas{t'}$ 
and $\ga\in\abs{t'}$. 

We can assume that $n\in\bas{a}$ and $\ga\in\abs{a}$, simply 
because $\sw ab$ is the identity outside of 
$\dom a=\bas{a}\ti\abs{a}$. 
We have $a^{t'}(n,\ga)=b(n,\ga)$ within this narrower domain, 
hence the result to prove is 
$q^{t'}_\ga\ogr n=h^{t}_{b(n,\ga)}$ for all $n\in\bas{a}$ 
and $\ga\in\abs{a}$. 
(Recall that $\sw ab$ does not change $h^t$, so that $h^{t'}=h^t$.)

However $q^{t'}_\ga\ogr n=q^{t}_\vt\ogr n$ by \ref{sQ1}, where 
$\vt=\swt abn(\ga)$, so that, in particular, $a(n,\vt)=b(n,\ga)$. 
Thus the equality required turns out to be 
$q^{t}_\vt\ogr n=h^{t}_{a(n,\vt)}$, which is true since $t$ is a 
condition.   

The other essential claim is that the action of $\sw ba$ is the 
inverse of the action of $\sw ab$. 
Suppose that $t\in\dT_a$ and let $t'=\sw ab\ap t$; $t\in\dT_b$.
Put $s=\sw ba\ap t'$; $s\in\dT_a$ once again. 
We have to show that $s=t$. 
The key fact is $q^s_\ga\og n=q^t_\ga\og n$ for all $n\in\bas a$ 
and $\ga\in\abs a$.
By definition $q^s_\ga\og n=q^{t'}_\za\og n$, where $\za=\swt ban$, 
in particular, $b(n,\za)=a(n,\ga)$. 
Still by definition, $q^{t'}_\za\ogr n=q^{t}_\vt\ogr n$, where 
$\vt=\swt abn(\za)$, so that $a(n,\vt)=b(n,\za)$.
To conclude, $q^s_\ga\og n=q^{t}_\vt\ogr n$, where
$a(n,\ga)=a(n,\vt)$.
But then $q^t_\ga\og n=q^{t}_\vt\ogr n$ by \ref{t2} of 
Section \ref{bd5}, and hence we have 
$q^s_\ga\og n=q^{t}_\ga\ogr n$, as required.

Claims \ref{sk2i}, \ref{sk2ii} are rather obvious.

It follows from \ref{sQ2} that 
claim \ref{sk2iii} is trivial for $n\in\abs t\bez B$, 
while in the case $n\in B$ it suffices to prove 
$\ens{q^{t'}_\ga\og n}{\ga\in B}=\ens{q^{t}_\ga\og n}{\ga\in B}$. 
The inclusion $\sq$ holds because 
$q^{t'}_\ga\og n=q^t_\vt\og n$ by \ref{sQ1}, where 
$\vt=\swt abn(\ga)$. 
The inclusion $\qs$ holds by the same reason with respect to 
the inverse swap $\sw ba$.  
\epf

\parf{Rotations}
\las{sym3}

This is a more complicated type of transformations, and we 
have to define it by extension beginning from most elementary 
conditions.

\punk{Simple rotations}
\las{sym3a}
 
If $d\in\dD$ and $p\in\dP$, 
or generally even $d\in\dda$ and $p\in\pP$, 
then define 
$d\ap p=p':\dom{p'}\to2$ so that $\dom p=\dom{p'}$ and
$$
p'(\al)=\left\{
\bay{rcl}
p(\al) &\text{whenever}& \al\in(\dom p)\bez d\,,\\[1ex]

1-p(\al)&\text{whenever}& \al\in d\cap\dom p\,.
\eay
\right. 
$$ 

Clearly $p\mapsto d\ap p$ is an order-preserving automorphism of 
$\dP$ and of $\pP$. 

Transformations of this type, as well as those based on them and 
defined below, will be called \rit{rotations}.

\punk{Rotations for narrow subconditions}
\las{sym3b}
 
We define product rotations which fit to 
conditions in $\pH$ and $\dH\sq\pH$. 
Let $\Psd$ consist of all indexed sets 
\imar{Psd}%
$\psi=\sis{\psi_\xi}{\xi\in\abs\psi}$, where
$\abs\psi\sq\omal$ is a finite set, and $\psi_\xi\in\ogd n$ 
for all $n\in\om$ and $\xi\in\abs\psi\cap\aal n$. 
If $\psi\in\Psd$ and $h\in\pH$ then define $h'=\psi\ap h\in\pH$ 
so that $\abs{h'}=\abs h$ and for all $\xi$:
$$
h'_\xi=\left\{
\bay{rcl}
h_\xi &\text{whenever}& \xi\in\abs h\bez \abs\psi\,,\\[1ex]

\psi_\xi\ap h_\xi &\text{whenever}& \xi\in \abs h\cap\abs\psi\,.
\eay
\right. 
$$ 
Let $\Psd\og n=\ens{\psi\in\Psd}{\abs\psi\sq\aal n}$; 
and accordingly if $\psi\in\Psi$ then let 
$\psi\og n=\psi\res\aal n$; then $\psi\og n\in\Psd\og n$.
The next lemma is obvious.

\ble
\lam{rl1}
If\/ $\psi\in\Psd$ then the map\/ $h\mapsto\psi\ap h$ is 
an order-preserving action\/ $\pH\onto\pH$ and\/ 
$\dH\onto\dH$.\qed
\ele

\punk{Rotations for wide subconditions}
\las{sym3c}
 
Now define product rotations which fit to conditions 
in $\pQ$ and $\dQ\sq\pQ$. 
Let $\Phd$ consist of all indexed sets 
\imar{Phd}%
$\vpi=\sis{\vpi_\xi}{\xi\in\abs\vpi}$, where
$\abs\vpi\sq\la$ is a finite set and $\vpi_\ga\in\dD$ 
for all $\ga\in\abs\vpi$. 
If $\vpi\in\Phd$ and $q\in\pQ$ then define $q'=\vpi\ap q\in\pQ$ 
so that $\abs{q'}=\abs q$ and for all $\ga$:
$$
q'_\ga=\left\{
\bay{rcl}
q_\ga &\text{whenever}& \ga\in\abs q\bez\abs\vpi\,,\\[1ex]

\vpi_\ga\ap q_\ga &\text{whenever}& \xi\in \abs\vpi\cap\abs q\,.
\eay
\right. 
$$
The lext elementary lemma is left to the reader. 

\ble
\lam{rl2}
If\/ $\vpi\in\Phd$ then the map\/ $q\mapsto\vpi\ap q$ is 
an order-preserving action\/ $\pQ\onto\pQ$ and\/ 
$\dQ\onto\dQ$. 
If\/ $q\in\pQ$ then\/ $q$ and\/ $\vpi\ap q$ are equally 
shaped.\qed
\ele

As above, say that $\vpi\in\Phd$ is 
\rit{compatible} with an assignment $a\in \dA$, 
in symbol $\vpi\in\Phd_a$, iff 
$\vpi_\ga\og n=\vpi_\da\og n$ holds 
whenever $\ga,\da\in\abs\vpi\cap\abs a$, $n\in\bas a$, and 
$a(n,\ga)=a(n,\da)$.
In this case, if in addition $\abs \vpi\sq\abs a$ 
then we define:
\ben
\aenu
\itla{ro1} 
a rotation $\psi=\ftp\vpi a\in\Psd$ 
(\dd a\rit{projection})
\imar{ftp vpi a}%
so that 
$$
\abs\psi=\ens{a(n,\ga)}{n\in\bas a\land\ga\in\abs\vpi}
$$ 
and if $n\in\bas a$, $\ga\in\abs\vpi$, and 
$\xi=a(n,\ga)$ then $\psi_\xi=\vpi_\ga\og n$;

\itla{ro2} 
a rotation $\ve=\ftw\vpi a\in\Phd$   
\imar{ftw vpi a}%
(\dd a\rit{extension})
so that $\abs{\ve}=\abs a$, 
$\ve_\da=\vpi_\da$ for all $\da\in\abs\vpi$, and the 
following holds for all $\ga\in\abs a\bez\abs \vpi$ and $n\in\om$:
\een
$$
\ve_\ga\og n=
\left\{
\bay{rcl}
\vpi_\da\og n &\text{iff}& 
n\in\bas a\,\land\,\da\in\abs\vpi\,\land\,a(n,\ga)=a(n,\da)\,,\\[1ex]
\pu &\text{iff}& 
n\nin\bas a\,\lor\,\neg\:
\sus\da\in\abs\vpi\,(a(n,\ga)=a(n,\da))\,.
\eay
\right.
$$
The consistency of both \ref{ro1} and \ref{ro2} 
follows from the compatibility assumption.

\punk{Rotations for conditions}
\las{sym3d}
 
Finally we define how any $\vpi\in\Phd$ acts
on the set  
$$
\dT_\vpi=\ens{t\in\dT}{\abs\vpi\sq\abs t\,\land\,
\vpi\text{ is compatible with }a^t}\,.
$$ 
If $t\in\dT_\vpi$ then let $\vpi\ap t=t'$, where 
$q^{t'}=(\ftw\vpi{a^t})\ap q^t$, 
$a^{t'}=a^t$, 
$h^{t'}=(\ftp\vpi{a^t})\ap h^t$.

\ble
\lam{sk3}
Suppose that\/ $\vpi\in\Phd$. 
Then the map\/ $t\mapsto\vpi\ap t$ is 
an order-preserving action\/ $\dT_\vpi\onto\dT_\vpi$, with\/ 
$t\mapsto\vpi^{-1}\ap t$ being the inverse.

If\/ $t\in\dT_\vpi$ is uniform then so is\/ $t'=\vpi\ap t$, 
and\/ $q^t,q^{t'}$ are equally shaped.
\ele
\bpf
Assume that $t\in\dT_\vpi$ and prove that $t'=\vpi\ap t$ 
belongs to $\dT_\vpi$ as well; this is the only part of the 
lemma not entirely trivial. 
We have to check \ref{t3} of Section~\ref{bd5}, that is, 
$h^{t'}_{a^{t'}(n,\ga)}=q^{t'}_\ga\ogr n$ 
for all $n\in\bas{t'}$ and $\ga\in \abs{t'}$.
By definition $a^{t'}=a^t$, $\bas{t'}=\bas{t}$, and 
$\abs{t'}=\abs{t}$, hence we have to prove   
$q^{t'}_\ga\ogr n=h^{t'}_{a^t(n,\ga)}$, 
for all $n\in\bas{t}=\bas{t'}$, $\ga\in\abs{t}=\abs{t'}$.

Note that $q^{t'}=(\ftw\vpi{a^t})\ap q^t$ and 
$h^{t'}=\psi\ap h^t$, where 
$\psi=\ftp\vpi{a^t}\in\Psd$.\vom

{\it Case 1\/}: $\ga\in\abs\vpi$. 
Then $q^{t'}_\ga\og n=\vpi_\ga\og n\ap q^t_\ga\og n$.
Let $\xi=a^t(n,\ga)$. 
By definition $h^{t'}_\xi=\psi_\xi\ap h^t_\xi$. 
On the other hand, $\psi_\xi=\vpi_\ga\og n$ and 
$h^t_\xi=q^t_\ga\og n$.
Therefore 
$h^{t'}_\xi=\vpi_\ga\og n\ap q^t_\ga\og n=q^{t'}_\ga\og n$, 
as required.\vom

{\it Case 2\/}: 
$\ga\nin\abs\vpi$, and there is an ordinal $\da\in\abs \vpi$ 
such that $a^t(n,\ga)=a^t(n,\da)$. 
Then the extended rotation $\ve=\ftw\vpi{a^t}$ satisfies 
$\ve_\ga\og n=\vpi_\da\og n$, and hence 
$q^{t'}_\ga\og n=\ve_\ga\og n\ap q^t_\ga\og n=
\vpi_\da\og n\ap q^t_\da\og n=q^{t'}_\da\og n=h^{t'}_\xi$, 
where $\xi=a^t(n,\ga)=a^t(n,\da)$ 
(we refer to Case 1), as required.\vom

{\it Case 3\/}: 
$\ga\nin\abs\vpi$, but there is no ordinal $\da\in\abs \vpi$ 
such that $a^t(n,\ga)=a^t(n,\da)$.
The extended rotation $\ve=\ftw\vpi{a^t}$ satisfies 
$\ve_\ga\og n=\pu$ in this case, and hence 
$q^{t'}_\ga\og n=q^t_\ga\og n$. 
Moreover, the Case~3 assumption means that 
$\xi=a^t(n,\ga)\nin \abs\psi$, and hence   
$h^{t'}_\xi=h^t_\xi$, and we are done.
\epf

\parf{The symmetry lemma}
\las{mst}

We begin with auxiliary definitions.
If $u\in\dT$ then let 
$$
\tle{u}=\ens{u'\in\dT}{u'\le u}\,.
$$

\bdf
\lam{coh2}
Suppose that $N\sq\om$ and\/ $\Ga\sq\la$ are finite sets. 
Conditions $s,t\in \dT$ are \rit{similar on\/ $N\ti\Ga$} iff 
\ben
\xaenu
\itlm{m1.}\msur
$\Ga\sq\abs s=\abs t$, $N\sq\bas s=\bas t$, 

\itlm{m2.}\msur
${q^s\res\Ga}={q^t\res\Ga}$ and the restricted assignments\/ 
$a^s\pes\Ga$ and\/ $a^t\pes\Ga$  are coherent 
(see Section~\ref{bd2}),

\itlm{m3.} 
if\/ $n\in N$ then\/ $h^s\og n=h^t\og n$, and\/  
$a^s(n,\ga)=a^t(n,\ga)$ for all\/ $\ga\in\Ga$,
\een
and \rit{strongly similar on\/ $N\ti\Ga$} if in addition
\ben
\xaenu
\atc3
\itlm{m1,}\msur
$s,t$ are uniform conditions, and 
$q^s,q^t$ are equally shaped (see Section~\ref{bd4}),

\itlm{m4}\msur
$\ran {a^s}=\ran {a^t}$ and $\abs{h^s}=\abs{h^t}$,

\itlm{m5} 
conditions $h^s$ and $h^t$ are regular at every $n\in\bas s\bez N$ 
(Section~\ref{bd3}),

\itlm{m6}\msur
$\ens{h^s_\xi}{\xi\in\abs{h^s}}=\ens{h^t_\xi}{\xi\in\abs{h^t}}$ 
--- then easily\/\\
$\ens{h^s_\xi}{\xi\in\abs{h^s}\cap\aal n}=
\ens{h^t_\xi}{\xi\in\abs{h^t}\cap\aal n}$ for all\/ $n$.\qed
\een
\eDf

\bte
[the symmetry lemma]
\lam{main}
Suppose that\/ $N\sq\om$, $\Ga\sq\la$ are finite sets, 
conditions\/ $s,t\in \dT$ are strongly similar on\/ $N\ti\Ga$, 
$B=\bas s=\bas t$, $\Da=\abs s=\abs t$. 
Then$:$
\ben
\renu
\itlm{mai1}
there exists a transformation\/ $\pi\in\pif$ such that\/ 
$\pi\og n$ is the identity for all\/ $n\in N$, 
condition\/ $u=\pi\ap s$ is strongly similar to\/ $t$ on\/ 
$N\ti\Ga$, and moreover\/ 
$\pi\ap h^s=h^u=h^t$, and\/ $a^u\pes\Ga=a^t\pes\Ga\,;$

\itlm{mai2}
condition\/ $v=\sw{a^u}{a^t}\ap u$ is strongly similar to\/ $t$ 
on\/ $N\ti\Ga$, and moreover\/ $h^v=h^u$ and\/ $a^v=a^t\,;$

\itlm{mai3}
there is a rotation\/ $\vpi\in\Phd_{a^v}$ 
{\rm(\ie, compatible with $a^v$)} such that\/ 
$\abs\vpi=\Da$, $\vpi_\ga\og n=\pu$ for all\/ 
$n\in B$ and\/ $\ga\in\Da$,\snos
{Then obviously $\vpi$ is compatible with each of the 
assignments $a^s,a^t,a^u,a^v$.}
and moreover\/ $t=\vpi\ap v\,;$ 

\itlm{mai4}\msur
$\tau= \vpi\circ \sw{a^u}{a^t}\circ \pi$  
is an order preserving bijection from\/ $\tle{s}$ 
onto\/ $\tle{t}\,;$ 

\itlm{mai5}
any condition\/ 
$s'\in\tle{s}$ is similar to\/ $t'=\tau\ap s'$ 
on\/ $N\ti \Ga$.
\een
\ete

\bpf
\ref{mai1} 
Let 
$\Xi=\abs{h^s}=\abs{h^t}$. 
Under our assumptions, obviously there is a 
transformation $\pi\in\pif$ such that  
\ben
\aenu
\itlm{pi1}\msur 
$\bas\pi=B$ and if $n\in N$ then $\pi\og n$ is the identity;

\itlm{pi2}\msur
$\pi(a^s(n,\ga))=a^t(n,\ga)$\snos
{As $s,t$ are similar on $\Ga$, here we avoid a contradiction 
related to the possibility of equalities 
$a^t(n,\ga)=a^t(n,\ga')$ for $\ga\ne\ga'$ in $\Ga$.} 
for all $n\in B$ and $\ga\in\Ga$;
\vyk{
\snos
{To have this property, it is important that 
\ref{m2} of Definition~\ref{coh2} holds for $p,q$.}
}

\itlm{pi3}\msur
$\pi$ maps the set $\Xi$ onto itself, 
and $\pi$ is the identity outside of $\Xi$,

\itlm{pi4}
if $\xi\in \Xi=\abs{h^s}$ then $h^s_\xi=h^t_{\pi(\xi)}$.
\een
The only point of contention is whether \ref{pi2} does not 
contradict to \ref{pi4}. 
That is, we have to check that 
$h^s_{a^s(n,\ga)}=h^t_{a^t(n,\ga)}$.
Note that $h^s_{a^s(n,\ga)}=q^s_\ga\og n$ and 
$h^t_{a^t(n,\ga)}=q^t_\ga\og n$ by \ref{t3} of 
Section~\ref{bdf}. 
On the other hand $q^s_\ga=q^t_\ga$ by \ref{m2.} 
of Definition~\ref{coh2}, as required. 

\ble
\lam{ml1}
The transformation\/ $\pi$ satisfies\/ \ref{mai1} of the 
theorem, and in addition if\/ $s'\in\tle s$ then\/ 
$s'$ is similar to\/ $u'=\pi\ap s'$ on\/ $N\ti\Ga$. 
\ele
\bpf[Lemma]
Prove that $h^u=\pi\ap h^s$ is equal to $h^t$. 
(This is a fragment of \ref{mai1}.) 
We have $\abs{h^u}=\ens{\pi(\xi)}{\xi\in\abs{h^s}}=\Xi$ 
by \ref{pi3}, and $\abs{h^t}=\Xi$ as well. 
Thus it remains to prove that $h^u_\eta=h^t_\eta$ for any 
$\eta=\pi(\xi)\in \Xi$, where $\xi\in \Xi$.
Yet by definition (Section~\ref{sym1}) $h^u_\eta=h^s_\xi$, 
and $h^t_\eta=h^s_\xi$ by \ref{pi4}.

The equality $a^u\pes\Ga=a^t\pes\Ga$ follows from \ref{pi2} 
since $a^u(n,\ga)=\pi(a^s(n,\ga))$. 

Prove that any $s'\in T\zd s'\le s$, is 
similar to $u'=\pi\ap s'$ on $N\ti\Ga$.

Item \ref{m1.} of Definition~\ref{coh2} holds 
for the pair of conditions $s',u'$ simply because the action 
of any $\pi\in\pif$ preserves $\abs\cdot$ and $\bas$. 

Prove \ref{m2.}. 
We have $q^{s'}=q^{u'}$ because the action of $\pi$ does not 
change $q^{s'}$ at all.
To show the coherence of $a^{s'}\pes\Ga$ and $a^{u'}\pes\Ga$
suppose that $\ga,\da\in\Ga$, $n\in\om$, and 
$a^{s'}(n,\ga)=a^{s'}(n,\da)$, and prove that 
$a^{u'}(n,\ga)=a^{u'}(n,\da)$. 
(The inverse implication can be checked pretty the same way.)

Suppose first that $n\in B$. 
Then $a^{s'}(n,\ga)=a^{s}(n,\ga)$ and 
$a^{s'}(n,\da)=a^{s}(n,\ga)$, therefore 
$a^{s}(n,\ga)=a^{s}(n,\da)$. 
It follows that $a^{t}(n,\ga)=a^{t}(n,\da)$
by the coherence in \ref{m2.} for $s,t$, therefore 
$a^{u}(n,\ga)=a^{u}(n,\da)$ since $a^u\pes\Ga=a^t\pes\Ga$, 
and finally $a^{u'}(n,\ga)=a^{u'}(n,\da)$, as required. 

Now suppose that $n\nin B$. 
Then the equality $a^{s'}(n,\ga)=a^{s'}(n,\da)$ implies 
$\ga=\da$ by \ref{a4} of Section~\ref{bdf}, so obviously  
$a^{u'}(n,\ga)=a^{u'}(n,\da)$.
 
To check \ref{m3.}, that is, $h^{u'}\og n=h^{s'}\og n$ 
and $a^{u'}(n,\ga)=a^{s'}(n,\ga)$ for all $\ga\in\Ga$ and $n\in N$, 
use the fact that $\pi\og n$ is the identity 
for any $n\in N$ by \ref{pi1}.

Prove that $s$ is strongly 
similar to $u=\pi\ap s$ on $N\ti\ga$.
We have \ref{m1,} of Definition~\ref{coh2} 
(for the pair of conditions $s',u'$) by rather obvious reasons. 
The equalities $\ran a^{u'}=\ran a^{s'}$ and 
$\abs{h^{u'}}=\abs{h^{s'}}$ in \ref{m4} hold by 
\ref{pi3} since $\ran a^{u'}$ is equal to the \dd\pi image 
of $\ran a^{s'}$. 
Finally the equality 
$\ens{h^{u'}_\xi}{\xi\in\abs{h^{u'}}}=
\ens{h^{s'}_\xi}{\xi\in\abs{h^{s'}}}$ 
in \ref{m6} holds whenever $u'=\pi\ap s'$ for some $\pi$.
We conclude that conditions $u$ and $t$ are 
strongly similar on $N\ti\Ga$. 
\epF{Lemma}


\ref{mai2}
Let $a=a^u$ and $b=a^t$. 
Thus $a,b\in\dA$, $\dom a=\dom b=B\ti \Da$,  
$\ran a=\ran b$, and $a\pes\Ga=b\pes\Ga$ by the above. 
Thus, as obviously $u\in \dtr_a$, we define 
$v=\sw ab\ap u\in \dtr_b$. 

\ble
\lam{ml2}
Condition\/ \ref{mai2} of the theorem holds, and in addition 
if\/ $u'\in\tle u$ then\/ $u'$ is similar to\/ 
$v'=\sw ab\ap u'$ on\/ $N\ti\Ga$. 
\ele
\bpf[Lemma]
That equalities $h^v=h^u$ and $a^v=a^t$ in \ref{mai2} hold 
is clear by definition: for instance swaps do not change $h^u$ 
at all. 

Prove that any $u'\in T\zd u'\le u$, is 
similar to $v'=\sw ab\ap u'$ on $N\ti\Ga$.

By definition (see Section~\ref{sym2}) $v'$ and $u'$ are 
equal outside of the domain $N\ti\Da$, and $h^{v'}=h^{u'}$. 
Therefore we can \noo\ assume that $\abs{v'}=\abs{u'}=\Da$ and 
$\bas{v'}=\bas{u'}=B$.
Then $a^{v'}=b=a^t$ and $a^{u'}=a=a^u$, thus the restricted 
assignments $a^{v'}\pes\Ga=b\pes\Ga$ and 
$a^{u'}\pes\Ga=a\pes\Ga$ are not merely coherent 
(as required by \ref{m2.} of Definition~\ref{coh2}) 
but just equal by the above. 
The equality $q^{v'}\res\Ga=q^{u'}\res\Ga$ in \ref{m2.} 
follows from $a\pes\Ga=b\pes\Ga$ as well.
And finally we have $h^{v'}=h^{u'}$ 
($\sw ab$ does not change this component), proving \ref{m3.}. 

Now prove that any $u$ is strongly 
similar to $v=\sw ab\ap u$ on $N\ti\ga$.
We skip \ref{m1,} of Definition~\ref{coh2} as clear and 
rather boring. 
Further, as $h^{v}=h^{u}$, we have $\abs{h^{v}}=\abs{h^{u}}$
in \ref{m4} and the whole of \ref{m6}. 
It remains to show $\ran a^v=\ran a^u$ in \ref{m4}. 
Recall that $a^v=a^t$ while conditions $s,t,u$ are strongly similar, 
therefore $\ran a^v=\ran a^t=\ran a^s=\ran a^u$. 
We conclude that conditions $v$ and $t$ are 
strongly similar on $N\ti\Ga$. 
\epF{Lemma}

\ref{mai3}
Thus $v,t$ are uniform conditions, strongly similar on $N\ti\Ga$, 
and $a^v=a^t$. 
In particular $q^v$ and $q^t$ are equally shaped, that is, 
in this case, $\abs{q^v}=\abs{q^t}=\Da$ and   
$\dom{q^v_\ga\og n}=\dom{q^t_\ga\og n}$ 
holds for all $\ga\in \Da$ and $n\in\om$.
Define a rotation $\vpi\in\Phd$ so that still $\abs\vpi=\Da$, 
and 
$$
\vpi_\ga\og n=\ens{\al\in\dom{q^v_\ga\og n}=\dom{q^t_\ga\og n}}
{q^v_\ga(\al)\ne q^t_\ga(\al)}
$$ 
for all $\ga\in \Da$ and $n\in\om$. 
Then clearly $\vpi\ap q^v=q^t$. 
Moreover $\vpi$ is compatible with $a^v=a^t$, because so are 
$q^t$ and $q^v$ in the sense of \ref{t2} of Section~\ref{bdf}. 
Thus conditions $v$ and $t$ belong to $\dT_\vpi$, so 
$\vpi\ap v$ makes sense. 

\ble
\lam{ml3}
Condition\/ \ref{mai3} of the theorem holds, and in addition 
if\/ $v'\in\tle v$ then\/ $v'$ is similar to\/ 
$t'=\vpi\ap v'$ on\/ $N\ti\Ga$.
\ele
\bpf[Lemma]
Recall that $a^v=a^t$ and $h^v=h^u=h^t$ by \ref{mai1}, \ref{mai2}. 
It follows by \ref{t3} of Section~\ref{bdf} that 
$q^v_\ga\og n=q^t_\ga\og n$, and hence $\vpi_\ga\og n=\pu$, 
whenever $\ga\in\Da$ and $n\in B$. 
To accomplish the proof of \ref{mai3} check that $\vpi\ap v=t$. 
Indeed $a^v=a^t$ since $\vpi$ does not change this component. 
Further, $q^t=\vpi\ap q^v$ simply by the choice of $\vpi$. 
Let us show that $h^t=h^v$ as well. 
Indeed, since by definition $\bas{h^v}=B=\bas t$, any change in 
$h^v$ by the action of $\vpi$ can be only due to a component 
$\vpi_\ga\og n$ for some $\ga\in\Da$ and $n\in B$ --- but this 
is the identity since $\vpi_\ga\og n=\pu$ in this case.

Now prove that any $v'\in T\zd v'\le v$, is 
similar to $t'=\vpi\ap v'$ on $N\ti\Ga$.
By definition $a^{v'}=a^{q'}$, covering the coherence in 
\ref{m2.} of Definition~\ref{coh2}. 
Further the extended rotation $\vpi'=\ftw\vpi{a^{v'}}$ 
obviously satisfies the same property 
$\vpi'_\ga\og n=\pu$ for all $n\in B$ and $\ga\in\Da'=\abs{v'}$.
This implies $h^{t'}\og n=h^{v'}\og n$ even for all $n\in B$, 
so that \ref{m3.} holds for $v',t'$ for all $n\in B$. 
It only remains 
to prove that $q^{t'}\res \Ga=q^{v'}\res\Ga$ in \ref{m2.} of 
Definition~\ref{coh2}, that is, $q^{t'}_\ga=q^{v'}_\ga$ 
for all $\ga\in\Ga$.  

By definition it suffices to show that $\vpi_\ga\og n=\pu$ 
for all $\ga\in\Ga$ and $n\in\om$, or equivalently, 
$q^v\res\Ga=q^t\res\Ga$ --- yet this is the case since 
$v$ and $t$ are similar on $N\ti\Ga$ by the above.
\epF{Lemma}

Finally, \ref{mai4} of the theorem 
is a consequence of lemmas \ref{sk1}, \ref{sk2}, \ref{sk3},
while \ref{mai5} is a corollary of lemmas \ref{ml1}, \ref{ml2},  
\ref{ml3}.\vtm 

\epF{Theorem}

\parf{The extension}
\las{E}

Let a set $G\sq\dT$ be \dd\dT generic over $\bL$.
It naturally produces:
\bit
\itsep
\item[--\ ]
for any $n$ and $\xi\in\aal n$, 
$\fg_\xi=\bigcup_{t\in G}h^t_\xi\in2^{\aal n}$,

\item[--\ ]
for every $n$, 
$\fg\og n=\sis{\fg_\xi}{\xi\in\aal n}:\aal n\to 2^{\aal n}$,

\item[--\ ]
for any $\ga<\la$, $\yg_\ga=\bigcup_{t\in G}q^t_\ga\in2^{\omal}$,

\item[--\ ]\msur
for any $\ga<\la$ and $n$, 
$\yg_\ga\og n=\yg_\ga\res\aal n\in2^{\aal n}$,

\item[--\ ]\msur
$\vy[G]=\sis{\yg_\ga}{\ga<\la}$, a map $\la\to 2^{\omal}$,

\item[--\ ]
a map $\ag=\bigcup_{t\in G}a^t:\om\ti\la\to\omal$ such that 
$\ag(n,\ga)\in\aal n$ for all $n$ and $\ga$. 
\eit

\ble
\lam{str}
If a set\/ $G\sq\dT$ is\/ \dd\dT generic over\/ $\bL$ then\/ 
\ben
\renu
\itsep
\itla{str1}
if\/ $n<\om$, $\ga<\la$, and\/ $\ag(n,\ga)=\xi$ 
then\/ $\yg_\ga\og n=\fg_\xi\;;$

\itla{str2}
if\/ $n<\om$, $\ga,\da<\la$, and\/ $\ag(n,\ga)\ne\ag(n,\da)$ 
then\/ $\yg_\ga\og n\ne\yg_\da\og n\;;$

\itla{str3}
if\/ $\ga\ne\da<\la$ then there is a number\/ 
$n_0=n_0(\ga,\da)$ such that\/
$\ag(n,\ga)\ne \ag(n,\da)$ for all\/ 
$n\ge n_0$. 
\een
\ele
\bpf
\ref{str1} is obvious.
 
\ref{str2}
Suppose that a condition $t\in G$ forces otherwise, and 
$\ga,\da\in\abs t$, $n\in\bas t$. 
Then $\xi=a^t(n,\ga)\ne a^t(n,\da)=\eta$; $\xi,\eta$ are 
ordinals in $\aal n$. 
Note that $h^t_\xi$ and $h^t_\eta$ are conditions in $\dP\og n$. 
Let $w_\xi\le h^t_\xi$ and $w_\eta\le h^t_\eta$ be any pair of 
{\ubf incompatible} conditions in $\dP\og n$. 
Let $t'\in T$ be a condition which differs from $t$ only in the 
following: $q^{t'}_\ga\og n=h^{t'}_\xi=w_\xi$ and 
$q^{t'}_\da\og n=h^{t'}_\eta=w_\eta$. 
Obviously $t'\le t$, and $t'$ forces that 
$\yg_\ga\og n\ne\yg_\da\og n$.

\ref{str3} 
Definitely there is a condition $t\in G$ such that $\abs t$ 
contains both $\ga$ and $\da$.
Let $B=\bas t$ (a finite subset of $\om$) 
and let $n_0$ be bigger than $\tmax B$. 
Now if $s\in G$, $s\le t$, and $n\in\bas s$, $n\ge n_0$, 
then $a^s\le a^t$, and hence $a^s(n,\ga)\ne a^s(n,\da)$.
This implies $\ag(n,\ga)\ne \ag(n,\da)$.
\epf

Now let us define a {\ubf symmetric subextension} of 
$\bL[G]$, on the base of certain 
symmetric hulls of sets $\fg\og n$ and $\yg_\ga$. 

\bbl
\lam{bla1}
Below, $\pif\zd\Phd\zd\Psd\zd\dD\zd\dD\og n\zd$ 
mean objects defined in $\bL$ 
as in Sections \ref{bd} --- \ref{sym3}. 
Thus in particular $\pif\in\bL$ and all elements of $\pif$ 
belong to $\bL$ either.\qed
\ebl

In $\bL[G]$, put
\bit
\item[--\ ]
for every $n$, 
$\Fg\og n=$ the \dd{(\pif,\Psd)}hull of $\fg\og n$. 
Thus the set $\Fg\og n$ consists of elements of the form 
$\pi\ap(\psi\ap \fg\og n)$, where $\pi\in\pif$ 
and $\psi\in\Psd$.

\item[--\ ]
$\vfg=\sis{\Fg\og n}{n<\om}$.
\eit
The actions of $\pi\in\pif$ and $\psi\in\Psd$ 
are defined as in sections \ref{sym1} and \ref{sym3} above. 
In particular $\psi\ap \fg\og n$ and $\pi\ap(\psi\ap \fg\og n)$ 
are maps $\aal n\to2^{\aal n}$ in $\bL[\fg\og n]$. 
It is clear that $\Fg\og n$ is closed under further 
application of transformations in $\pif$ and $\Psd$, 
so there is no need to consider iterated actions. 

It takes more time to define suitable hulls of 
elements $\yg_\ga$. 
First of all, put
\bit
\itsep
\item[--\ ]
for any $n$ and $\ga<\la$, 
$\Yg_\ga\og n=\ens{d\ap \yg_\ga\og n}{d\in\dD\og n}\sq 
2^{\aal n}$; 

\item[--\ ]
for any $n$,   
$\Yg\og n=\bigcup_{\ga<\la}\Yg_\ga\og n$ --- 
still $\Yg\og n\sq 2^{\aal n}$, and obviously 
$\Yg\og n$ is the \dd{\dD\og n}hull of 
$\ens{\yg_\ga\og n}{\ga<\la}$.
\eit
Finally, if $\ga<\la$ then we let $\Yg_\ga$ be the set of all 
$z\in 2^{\omal}$ in $\bL[G]$ such that there exist a set 
$d\in\dD$ 
and a number $n_0$ satisfying:\vtm

1) $z\og n=d\og n\ap \yg_\ga\og n$ for all $n\ge n_0$;\snos
{Regarding the action of $d\in\dD$ see Section~\ref{sym3a}.}\vtm

2) $z\og n\in \Yg\og n$ for all $n<n_0$.\vtm

\noi
In other words, to obtain $\Yg_\ga$ we first define the 
\dd\dD hull $\dD\ap{\yg_\ga}=\ens{d\ap\yg_\ga}{d\in D}$ of 
$\yg_\ga$, and then allow to substitute sets in $\Yg\og n$ 
for $y\og n$ for any $y\in \dD\ap{\yg_\ga}$ and finitely 
many $n$, so that
\bit
\item[$(\star)$]\msur
$\Yg_\ga$ is the set of all $z\in 2^{\omal}$ (in $\bL[G]$) such 
that there exist an element $y\in\dD\ap{\yg_\ga}$ 
and a number $n_0$ satisfying: 
$z\og n=y\og n$ for all $n\ge n_0$, and 
$z\og n\in \Yg\og n$ for all $n<n_0$.
\eit

\ble
\lam{cons}
If\/ $\ga\ne\da$ then\/ $\Yg_\ga\cap\Yg_\da=\pu$.
\ele
\bpf
Suppose towards the contrary that $z\in\Yg_\ga\cap\Yg_\da$.
Then by $(\star)$ there exist rotations $d',d''\in\dD$ and a 
number $n_0$ such that the elements $y'=d'\ap\yg_\ga$ and 
$y''=d''\ap\yg_\da$ satisfy $y\og n=y'\og n$ for all $n\ge n_0$. 
In other words, $\yg_\ga\og n=(d\ap\yg_\da)\og n$ for all 
$n\ge n_0$, where $d=d'\sd d''\in\dD$ (symmetric difference).
Now use Lemma~\ref{str}\ref{str3} to obtain a number $n\ge n_0$ 
such that $\ag(n,\ga)\ne\ag(n,\da)$; still we have 
$\yg_\ga\og n=d\og n\ap\yg_\da\og n$.
But this yields a contradiction similarly to the proof of 
Lemma~\ref{str}\ref{str2}.
\epf

Now we let, in $\bL[G]$, 
$\vyg=\sis{\Yg_\ga}{\ga<\la}$, a function  
defined on $\la$.

We finally define
$$
\textstyle
\wg=\bigcup_{n}\Fg\og n\,\cup\,
\bigcup_{\ga<\la}\Yg_\ga
\,\cup\,\ans{\vfg,\vyg}\,.
$$

\bdf
\lam{se}
$\Ls[G]=\bL(\wg)=$ HOD over $\wg$ in $\bL[G]$.
\edf

Thus by definition every set in $\Ls G$ is definable in 
$\bL[G]$ by a formula with parameters in $\bL$, 
two special parameters $\vfg$ and $\vyg$, 
and finally parameters which belong to the sets 
$\Fg\og n$ and $\Yg_\ga$ for various $n<\om$ and $\ga<\la$. 
The next lemma allows to reduce the last category of 
parameters, basically, to those in 
$\ens{\fg\og n}{n<\om}\cup\ens{\yg_\ga}{\ga<\la}$.

\ble
\lam{redp}
If\/ $n<\om$ then every\/ $x\in\Fg\og n$ belongs to\/ 
$\bL[\fg\og n]$. 
If\/ $\ga<\la$ and\/ $z\in\Yg_\ga$ then there is a finite 
set\/ $\Da\sq\la$ such that\/ $z\in\bL[\ens{\yg_\da}{\da\in\Da}]$.
\ele
\bpf
By definition $x$ belongs to the \dd{(\pif,\Psd)}hull of 
$\fg\og n$. 
But $\pif$ and $\Psd$ belong to $\bL$ 
(see Blanket Agreement~\ref{bla1}). 
Regarding the claim for $z\in\Yg_\ga$, come back to $(\star)$. 
Note that $y$ as in $(\star)$ belongs to $\bL[\yg_\ga]$ 
(since $\dD\in\bL$). 
Then to obtain $z$ from $y$ we replace a finite number of 
intervals $y\og n$ in $y$ by elements of sets $\Yg\og n$. 
Thus suppose that $n<\om$ and $w\in\Yg\og n$, that is, 
$w\in\Yg_\da\og n$, where $\da<\la$. 
But then $w\in\bL[\yg_\da]$ (since $\dD\og n\in\bL$), so that 
it suffices to define $\Da$ as the (finite) set of all ordinals 
$\da$ which appear in this argument for all intervals $y\og n$ 
to be replaced. 
\epf

\parf{Definability lemma}
\las{pdl}

The next theorem plays key role in the analysis of the 
abovedefined symmetric subextension.

\bte
[the definability lemma]
\lam{dl1}
Suppose that a set\/ $G\sq\dT$ is\/ \dd\dT generic over\/ $\bL$, 
and\/ $N\sq\om\zt \Ga\sq\la$ are finite sets. 
Let\/ $Z\in\bL[G]\zt Z\sq\bL$, be a set definable in\/ $\bL[G]$ 
by a formula with parameters in\/ $\bL$ and those in the list\/ 
$$
\ans{\vfg,\vyg}\,\cup\,\ens{\fg\og n}{n\in N}
\,\cup\,\ens{\yg_\ga}{\ga\in\Ga}.
$$
Then\/ 
$Z\in\bL[\ens{\fg\og n}{n\in N},\ens{\yg_\ga}{\ga\in\Ga}]$.
\ete

Beginning the {\ubf proof of Theorem~\ref{dl1}}, we
put $\vxng=\sis{\fg\og n}{n\in N}$ and 
$\vyng=\sis{\yg_\ga}{\ga\in\Ga}$,
and let 
$$
\vt(z):=
\vt(z,\vfg,\vyg,\vxng,\vyng) 
$$ 
be a formula such that $Z=\ens{z}{\vt(z)}$ in $\bL[G]$.
%
By Lemma~\ref{str}\ref{str3} there is $n_0$ such that 
$\ag(n,\ga)\ne \ag(n,\da)$ whenever $n>n_0$ and $\ga\ne\da$ 
belong to $\Ga$. 

Let $M=N\cup\ens{n}{n\le n_0}$. 
Say that a condition $t\in\dT$ 
\rit{complies with\/ $\vxng\zi\vyng$}
if $M\sq\bas t$, $\Ga\sq\abs t$, and 
\ben
\itsep
\Renu
\itlm{pa1} 
if $n\in N$ and $\xi\in\abs{h^t}\cap\aal n$ then 
$h^t_\xi\subset \fg_\xi$,
 
\itlm{pa2} 
if $\ga\in\Ga$ then 
$q^t_\ga\subset \yg_\ga$,
 
\itlm{pa3} 
if $n\in\bas t$ and $\ga\in\Ga$ then
$a^t(n,\ga)=\ag(n,\ga)$.
\een
For instance any condition $t\in G$ with $M\sq\bas t$, 
$\Ga\sq\abs t$ complies with $\vxng\zi\vyng$ by obvious reasons. 

It is quite clear that 
the set $\txy$ of all conditions $t\in\dT$ which comply  
with $\vxng\zi\vyng$ belongs to $\bL[\vxng\zi\vyng]$. 
Therefore to prove the theorem it suffices to 
verify the following assertion:\vtm

\rit{if\/ $z\in\bL$, $s,t\in\txy$, and\/ $s$ forces\/ $\vt(z)$, 
then\/ $t$ does not force\/ $\neg\:\vt(z)$.}\vtm

\noi 
Suppose towards the contrary that this fails, so that
\ben
\fenu
\itlm*\msur
$z\in\bL$, $s,t\in\txy$, condition $s$ forces\/ $\vt(z)$, 
while $t$ forces $\neg\:\vt(z)$. 
\een

The proof of Theorem~\ref{dl1} continues in Sections \ref{pdl1} 
and \ref{pdl2}.

\parf{Proof of the definability lemma, part 1}
\las{pdl1}

Working towards the symmetry lemma.
Our goal is now to strengthen $s,t$ towards the requirements of 
Theorem~\ref{main}. 

\ble
\lam{dl2}
There exists a condition\/ $s'\in\txy$ such that\/ 
$\abs{s'}=\abs{s}\cup\abs t$, $\bas{s'}=\bas{s}\cup\bas t$,
and\/ $s'\le s$.
Accordingly there is a condition\/ $t'\in\txy$ such that\/ 
$\abs{t'}=\abs{s}\cup\abs t$, $\bas{t'}=\bas{s}\cup\bas t$, 
and\/ $t'\le t$.
\ele
\bpf[Lemma]
We define $a^{s'}$. 
This takes some time.\vom 

\rit{Domain\/ $\bas s\ti \abs s$.} 
If $n\in\bas s$ and $\ga\in\abs s$ then put 
$a^{s'}(n,\ga)=a^s(n,\ga)$ and  
$q^{s'}(n,\ga)=q^s(n,\ga)$.\vom 

\rit{Domain\/ $(\bas t\bez\bas s)\ti\Ga$.} 
If $n\in\bas t\bez\bas s$ and $\ga\in\Ga$ then put 
$a^{s'}(n,\ga)=\ag(n,\ga)$, and 
$q^{s'}(n,\ga)=q^s(n,\ga)$, as above.\vom 

\rit{Domain\/ $(\bas t\bez\bas s)\ti(\abs s\bez\Ga)$.} 
For any $n\in\bas t\bez\bas s$ fix a bijection $\da\mto\xi^n_\da$ 
from $\abs s\bez\Ga$ to $\aal n\bez \ens{a^{s'}(n,\ga)}{\ga\in\Ga}$.  
If now $\da\in\abs s\bez\Ga$ then put 
$a^{s'}(n,\da)=\xi^n_\da$ and $q^{s'}(n,\da)=\pu$.\vom 

\rit{Domain\/ $(\bas t\cup\bas s)\ti(\abs t\bez\abs s)$.} 
Fix an ordinal $\da^\ast\in\abs s$. 
If $n\in\bas t\cup\bas s$ and $\da\in\abs t\bez\abs s$ then put 
$a^{s'}(n,\da)=a^{s'}(n,\da^\ast)$ and 
$q^{s'}(n,\da)=q^{s'}(n,\da^\ast)$.\vom 

\rit{Domain\/ 
$\big(\om\bez{(\bas t\cup\bas s)}\big)\ti(\abs t\bez\abs s)$.} 
If $n\nin\bas t\cup\bas s$ and $\da\in\abs t\cup\abs s$ then put 
$q^{s'}(n,\da)=\pu$ and keep $a^{s'}(n,\da)$ undefined.\vom

On the top of the above definition, define $h^{s'}$ so that 
$$
\abs{h^{s'}}=\abs{h^{s}}
\cup\ens{\xi^n_\da}{n\in\bas t\bez\bas s\land \da\in\abs s\bez\Ga},
$$
$h^{s'}_\xi= h^{s}_\xi$ for all $\xi\in\abs{h^{s}}$, and 
$h^{s'}_{\xi^n_\da}=\pu$ for all
$n\in\bas t\cup\bas s$ and $\da\in\abs t\bez\abs s$.

We claim that $s'$ is as required. 

The key issue is to prove $a^{s'}\le a^s$, in particular, 
\ref{a4} of Section~\ref{bdf} for $a=a^{s'}$, $b=a^s$. 
Note that if $\ga\ne\da$ belong to $\Ga$ and $n\nin\bas s$ 
then $\ag(n,\ga)\ne \ag(n,\da)$ by the choice of $M$ and because 
$M\sq\bas s$. 
Therefore if $n\in\bas t\bez \bas s$ and $\ga,\da$ as indicated 
then by definition $a^{s'}(n,\ga)\ne a^{s'}(n,\da)$, as required. 

We have \ref{pa1}, \ref{pa2}, \ref{pa3} by obvious 
reasons: 
in particular, $q^{s'}_\ga=q^s_\ga$ for all $\ga\in\Ga$, 
and if $n\in N$ then $n\in\abs s$ and hence by construction 
$\abs{h^{s'}}\cap\aal n= \abs{h^{s}}\cap\aal n$ and 
$h^{s'}_\xi= h^{s}_\xi$ for all 
$\xi\in\abs{h^{s'}}\cap\aal n$.
\epF{Lemma}

It follows from the lemma that we can \noo\ assume in \ref* that 
\ben
\aenu
\itlm{st1}
conditions $s,t$ satisfy $\abs s=\abs t$ and $\bas s=\bas t$. 
\een

Moreover we can \noo\ assume that in addition to \ref* 
and \ref{st1}:
\ben
\atc1
\aenu
\itlm{st2}\msur
$\abs{h^s}=\abs{h^t}$, and 
if $n\in\bas s=\bas t$ then the set
$\abs{h^s}\cap\aal n=\abs{h^t}\cap\aal n$  
is infinite.
\een
This is rather elementary. 
If say $\xi\in\abs{h^s}\bez\abs{h^t}$ then simply add $\xi$ to 
$\abs{h^t}$ and define $h^t_\xi=\pu$.

Further, we can \noo\ assume that, in addition to \ref*,
\ref{st1}, \ref{st2}: 
\ben
\atc2
\aenu
\itlm{st4}
conditions $s,t$ satisfy $\ran a^s=\ran a^t$.    
\een
Suppose that $n\in\bas s$ and, say, 
$\xi\in (\ran a^t\bez \ran a^s)\cap\aal n$. 
Put $\xi_n=\xi$ and for any $m\in\bas s\zt m\ne n$ pick an 
ordinal $\xi_m\in\abs{h^s}\cap\aal m$, 
$\xi_m\nin\ran{a^s}\cup\ran{a^t}$ 
(this is possible by \ref{st2}). 
Add an ordinal $\ga\nin\abs s=\abs t$ to $\abs s$ and to 
$\abs t$. 
If $m\in\bas s=\bas t$ then put 
$a^s(m,\ga)=a^t(m,\ga)=\xi_m$ and   
$q^s_\ga\og m=q^t_\ga\og m=h^t_{\xi_m}$, 
and in addition define $q^s_\ga\og m=q^t_\ga\og m=\pu$ 
for all $m\nin\bas s=\bas t$. 
Conditions $s,t$ extended this way still satisfy 
\ref*, \ref{st1}, \ref{st2}, but now $\xi\in\ran{a^s}$. 
One has to maintain such extension for all indices 
$\xi$ in $\ran a^t\bez \ran a^s$ and 
$\ran a^s\bez \ran a^t$ one by one; 
the details are left to the reader.

\bre
After this step, the sets $\Da=\abs s=\abs t$ and $B=\bas s=\bas t$ 
(finite subsets of resp.\ $\la$ and $\om$) 
will not be changed, as well as the assignments 
$a=a^s$ and $b=a^t$ ($\dom a=\dom b=B\ti\Da$).
Put $\Xi=\abs{h^s}=\abs{h^t}$.
\ere

Further we can \noo\ assume that in addition to \ref*,
\ref{st1}, \ref{st2}, \ref{st4}:
\ben
\atc3
\aenu
\itlm{st5}
subconditions $q^s,q^t$ are uniform and equally 
shaped.    
\een
It suffices to define a pair of stronger conditions $s',t'\in\txy$ 
such that 
$$
\abs{s'}=\abs{t'}=\Da\,,\;
\abs{h^{s'}}=\abs{h^{t'}}=\Xi\,,\;
\bas{s'}=\bas{t'}=B\,,\; 
a^{s'}=a\,,\; 
a^{t'}=b\,,
$$ 
and in addition $q^{s'},q^{t'}$ are uniform 
and equally shaped.

Consider any $n\in\om$.
Put 
$d\og n=\bigcup_{\da\in\Da}(\dom q^s_\da\og n\cup\dom q^t_\da\og n)$, 
a set in $\dD\og n$.
If $\da\in\Da$ then define extensions 
$q^{s'}_\da\og n\zi q^{t'}_\da\og n\in\dP\og n$ of resp.\ 
$q^{s}_\da\og n\zi q^{t}_\da\og n$ so that
\ben
\renu
\itlm{zz1}\msur
$\dom{q^{s'}_\da\og n}=\dom{q^{t'}_\da\og n}=d\og n$, 

\itlm{zz2} 
if $n\in N$ and $\da\in\Ga$ then simply 
$q^{s'}_\da\og n= q^{t'}_\da\og n=\yg_\da\res{d\og n}$,

\itlm{zz3} 
if $n\in B$ and $\ga,\da\in\Da$ then: 
if $a(n,\da)=a(n,\ga)$ then $q^{s'}_\da\og n= q^{s'}_\ga\og n$,\\ 
and 
if $b(n,\da)=b(n,\ga)$ then $q^{t'}_\da\og n= q^{t'}_\ga\og n$.
\een
On the top of this, define $h^{s'}_{a(n,\da)}=q^{s'}_\da\og n$ 
and $h^{t'}_{b(n,\da)}=q^{t'}_\da\og n$ for all $n\in B$ and 
$\da\in \Da$. 
In the rest, put $\abs{h^{s'}}=\abs{h^{t'}}=\Xi$ 
(recall that $\Xi=\abs{h^s}=\abs{h^t}$), and 
$h^{s'}_\xi=h^s_\xi$, $h^{t'}_\xi=h^t_\xi$ 
for all $\xi\in\Xi$ {\ubf not} in 
$\ran a=\ran b$.

Further, we can \noo\ assume that, in addition to \ref* and 
\ref{st1} --- \ref{st5}: 
\ben
\atc4
\aenu
\itlm{st3}
conditions $s,t$ coincide on the domain $N\ti\Ga$, so that  
\ben
\itla{st3a} if $\ga\in\Ga$ then $q^s_\ga=q^t_\ga$,  

\itla{st3b} if $n\in N$ then $h^s\og n=h^t\og n$, that is, 
$h^s_\xi=h^t_\xi$ for all 
$\xi\in\abs{h^s}\cap\aal n=\abs{h^t}\cap\aal n$, and 

\itla{st3c} if $n\in N$ and $\ga\in\Ga$ then 
$a^s(n,\ga)=a^t(n,\ga)=\ag(n,\ga)$ --- but this already
follows from the compliance assumption. 
\een
\een
Regarding \ref{st3a}, note that this is already done. 
Indeed, $q^s,q^t$ are equally shaped by \ref{st5}, and 
satisfy $q^s_\ga\su \yg_\ga$ and $q^t_\ga\su \yg_\ga$ 
by \ref*, therefore $q^s_\ga=q^t_\ga$. 

Now consider \ref{st3b}; suppose that $n\in N$. 
Let $\xi\in\abs{h^s}\cap\aal n$. 

If $\xi\in\ran{a^s}=\ran{a^t}$ then $\xi=a^s(n,\ga)=a^t(n,\da)$
for some $\ga,\da\in\Da$, and then $h^s_\xi=q^s_\ga\og n$ and 
$h^t_\xi=q^t_\ga\og n$. 
It follows that $\dom{h^s_\xi}=\dom{h^t_\xi}$, by \ref{st5}. 
Therefore ${h^s_\xi}={h^t_\xi}$, because we have 
$h^s_\xi\su \fg_\xi$ and $h^t_\xi\su \fg_\xi$. 

If $\xi\in\ran{a^s}=\ran{a^t}$ then still 
$h^s_\xi\su \fg_\xi$ and $h^t_\xi\su \fg_\xi$, thus 
$h^s_\xi$ and $h^t_\xi$ are compatible as conditions in $\dP$, 
and we simply replace either of them by $h^s_\xi\cup h^t_\xi$.

And finally, we can \noo\ assume that, in addition to \ref* and 
\ref{st1} --- \ref{st3}: 

\ben
\atc5
\aenu
\itlm{st6}
we have    
$\ens{h^s_\xi}{\xi\in\abs{h^s}}=\ens{h^t_\xi}{\xi\in\abs{h^t}}$ 
as in \ref{m6} of Definition~\ref{coh2}, and 
subconditions $h^s,h^t$ are regular on every $n\in B\bez N$ 
(Subsection~\ref{bd3}).    
\een
The equality 
$\ens{h^s_\xi}{\xi\in\abs{h^s}\cap\aal n}=
\ens{h^t_\xi}{\xi\in\abs{h^t}\cap\aal n}$
holds already for all $n\in N$ by \ref{st3}. 

Now suppose that $n\in B\bez N$. 
The requirement of compliance with $\vxng\zi\vyng$ is void for 
$n\nin N$, therefore we can simply extend 
$h^s\og n$ and $h^t\og n$ to a 
bigger domain and appropriately define $h^s_\xi$ and $h^t_\xi$ 
for all ``new'' elements $\xi$ in these extended domains 
so that \ref{st6} holds, without changing $q^s\zi q^t$ and 
$a^s\zi a^t$.

To conclude, we can \noo\ assume in \ref* that 
\ref{st1} --- \ref{st6} hold, 
that is, in other words, conditions $s\zi t\in\txy$ are 
strongly similar on $N\ti\Ga$ in the sense of Definition~\ref{coh2}.

\parf{Proof of the definability lemma, part 2}
\las{pdl2}

We continue the proof of Theorem~\ref{dl1}. 
Our intermediate result and the starting point of the final 
part of the proof is the contrary assumption \ref* with 
the additional assumption that conditions $s\zi t\in\txy$ 
in \ref* are strongly similar on $N\ti\Ga$, 
and to complete the proof of the theorem it suffices to derive a 
contradiction. 
This will be obtained by means of Theorem~\ref{main}.

In accordance with Theorem~\ref{main}, let 
$B=\bas s=\bas t$, $\Da=\abs s=\abs t$, and let 
transformations $\pi\zi \sw{a^u}{a^t}\zi\vpi$ and 
$\tau=\vpi\circ \sw{a^u}{a^t} \circ \pi$, and conditions 
$v,u\in\dT$ satisfy $\bas u=\bas v=B$, $\abs u=\abs v=\Da$, and
\ben
\renu
\itlm{nai1}\msur
$\pi\in\pif$, 
$\pi\og n$ is the identity for all\/ $n\in N$, 
$u=\pi\ap s$, $u$ is strongly similar to\/ $t$ on\/ 
$N\ti\Ga$, and moreover\/ 
$\pi\ap h^s=h^u=h^t$, and\/ $a^u\pes\Ga=a^t\pes\Ga\,;$

\itlm{nai2}\msur
$v=\sw{a^u}{a^t}\ap u$, $v$ is strongly similar to\/ $t$ 
on\/ $N\ti\Ga$, $h^v=h^u$, $a^v=a^t\,;$

\itlm{nai3}\msur
$\vpi\in\Phd_{a^v}$, 
$\abs\vpi=\Da$, $\vpi_\ga\og n=\pu$ for all\/ 
$n\in B$ and\/ $\ga\in\Da$,
and\/ $t=\vpi\ap v\,;$ 

\itlm{nai4}\msur
$\tau= \vpi\circ \sw{a^u}{a^t}\circ \pi$  
is an order preserving bijection from\/ $\tle{s}$ 
onto\/ $\tle{t}\,;$ 

\itlm{nai5}
any condition\/ $s'\in\tle{s}$ is similar to\/ $t'=\tau\ap s'$ 
on\/ $N\ti \Ga$.
\een
(= items \ref{mai1} --- \ref{mai5} of Theorem~\ref{main}).

Consider a set $G\sq\dT$ generic over $\bL$ and containing $s$. 
We assume that $s$ is the largest (= weakest) condition in $G$. 
Then, by \ref{nai5}, $H=\ens{\tau\ap s'}{s'\in G}\sq\dT$ is 
generic over $\bL$ either, and $\bL[H]=\bL[G]$. 
Moreover $t=\tau\ap s\in H$.
Therefore it follows from \ref* that 
\ben
\fenu
\atc1
\itlm+\msur
$\vt(z,\vfg,\vyg,\vxng,\vyng)$ is true in $\bL[G]$, but\\[1ex]
$\vt(z,\vfh,\vyh,\vxnh,\vynh)$ is false in $\bL[H]=\bL[G]$.
\een
Our strategy to derive a contradiction will be to show that the 
parameters in the formulas are pairwise equal, and hence one 
and the same formula is simultaneously true and false in one and 
the same class. 
This is the content of the following lemma.

\ble
\label{peq}
\ben
\renu
\itlm{peq1}\msur 
$\vyng=\vynh$, that is, if\/ $\ga\in\Ga$
then\/ $\yg_\ga=\yh_\ga\;;$

\itlm{peq2}\msur 
$\vxng=\vxnh$, that is, if\/ $n\in N$ and\/ $\xi\in\aal n$
then\/ $\fg_\xi=\fh_\xi\;;$ 

\itlm{peq3}\msur 
$\vfg=\vfh$, that is, $\Fg\og n=\Fh\og n$ for all\/ $n\in\om\;;$ 

\itlm{peq4}\msur 
$\vyg=\vyh$, that is, $\Yg_\ga=\Yh_\ga$ for all\/ $\ga<\la$. 
\een
\ele
\bpf
\ref{peq1}
If $\ga\in\Ga$ then by definition  
$\yg_\ga=\bigcup_{s'\in G} q^{s'}_\ga$ and 
$\yh_\ga=\bigcup_{t'\in H} q^{t'}_\ga=
\bigcup_{s'\in G} q^{(\tau\ap s')}_\ga$. 
Yet if $s'\in G$ then condition $t'=\tau\ap s'$ satisfies 
$q^{t'}_\ga=q^{s'}_\ga$ by \ref{nai5}.

\ref{peq2}
A similar argument.
Suppose that $n\in N$ and $\xi\in\aal n$.
By definition, $\fg_\xi=\bigcup_{s'\in G}h^{s'}_\xi$ and 
$\fh_\xi=\bigcup_{t'\in H}h^{t'}_\xi=
\bigcup_{s'\in G}h^{(\tau\ap s')}_\xi$. 
However if $s'\in G$ then condition $t'=\tau\ap s'$ satisfies 
$h^{t'}_\xi=h^{s'}_\xi$ still by \ref{nai5}.

\ref{peq3}
By definition, $\Fg\og n$ and $\Fh\og n$ are the 
\dd{(\pif,\Psd)}hulls of resp.\ 
$$
\fg\og n=\sis{\fg_\xi}{\xi\in\aal n}
\quad\text{and}\quad
\fh\og n=\sis{\fh_\xi}{\xi\in\aal n}\,.
$$ 
Thus it remains to prove that $\fh\og n$ belongs to the 
\dd{(\pif,\Psd)}hull of $\fg\og n$, and vice versa. 

Let $\psi=\ftp\vpi{a^v}$ 
(a rotation in $\Psd$, see Section~\ref{sym3}). 
By definition, if $s'\in G$ and $t'=\tau\ap s'$, then the 
subconditions $h^{s'}$ and $h^{t'}$ satisfy 
$h^{t'}=\psi\ap(\pi\ap h^{s'})$. 
(The middle transformation $\sw{a^u}{a^t}$ does not act on the 
\dd hcomponents). 
It easily follows that $\fh\og n=\psi\ap(\pi\ap\fg\og n)$, as 
required.  

\ref{peq4}
Note that the sequences $\vyng=\sis{\yg_\ga}{\ga<\la}$ 
and $\vynh=\sis{\yh_\ga}{\ga<\la}$ satisfy 
$\vynh=\vpi\ap(\sw{a^u}{a^t}\ap\vyng)$.
(Permutation $\pi$ does not act on wide subconditions.) 
That is, the construction of $\vynh$ from $\vyng)$ goes in 
two steps.\vom

{\it Step 1\/}: 
we define $\vz=\sis{\zzz_\ga}{\ga<\la}$ by 
$\vz=\sw{a^u}{a^t}\ap\vyng$. 
Thus by definition\vtm

1) $\zzz_\ga\in2^{\omal}$ for all $\ga$,\vtm

2) 
if $n\nin B$ or $\ga\nin\Da$ then $\zzz_\ga\og n=\yg_\ga\og n$, 
and\vtm

3) 
\psur 
if $n\in B$ and $\ga\in\Da$ then $\zzz_\ga\og n=\yg_{\vt}\og n$, 
where $\vt=\swt{a^u}{a^t}n(\ga)$.\vtm

\noi
Thus the difference between $\vz$ and $\vyng$ is located within 
the finite domain $B\ti\Da$. 
Moreover, as in Lemma~\ref{sk2}\ref{sk2iii}, we have 
\ben
\fenu
\atc2
\itla{ddag}
$\ens{\zzz_\ga\og n}{\ga<\la}=\ens{\yg_\ga\og n}{\ga<\la}$ for 
every $n$.
\een

{\it Step 2\/}: 
we define $\vynh=\vpi\ap\vz$.   
Thus by definition\vtm

4) 
if $\ga\in\Da$ then directly $\yh_\ga=\vpi\ap\zzz_\ga$, 
that is, $\yh_\ga\og n=\vpi\og n\ap\zzz_\ga\og n$, $\kaz n$;\vtm

5) 
if $\ga\nin\Da$ and 
$\sus\da\in \Da\,(\ag(n,\ga)=\ag(n,\da))$, then 
$\yh_\ga\og n=\vpi\og n\ap\zzz_\ga\og n$;\vtm

6) 
if $\ga\nin\Da$ but   
$\neg\:\sus\da\in \Da\,(\ag(n,\ga)=\ag(n,\da))$, then 
$\yh_\ga\og n=\zzz_\ga\og n$.\vtm

\noi
Now it immediately follows from \ref{ddag} that $\Yg\og n=\Yh\og n$ 
for every $n$: both sets are equal to the \dd{\dD\og n}hull of one 
and the same set mentioned in \ref{ddag}.

We are ready to prove that $\Yg_\ga=\Yh_\ga$ for every 
$\ga<\la$.

We start with a couple of definitions. 
If $y,y'\in 2^{\aal n}$ and there exists a set 
$d\in\Yg\og n=\Yh\og n$ such that $y'=d\ap y$ then write 
$y\eqn n y'$.
If $y,y'\in 2^{\omal}$ and there exists a number $n_0$  
such that $y'\og n\eqn n y\og n$ for all $n<n_0$ and 
$y'\og n=y\og n$ for all $n\ge n_0$ then write 
$y\eqa y'$.

Then by $(\star)$ in Section~\ref{E} we have:
$$
\left.
\bay{rcl}
\Yg_\ga &=&
\ens{z\in2^{\omal}}{\sus y\in\dD\ap\yg_\ga\,(z\eqa y)};\\[1ex]
\Yh_\ga &=&
\ens{z\in2^{\omal}}{\sus y\in\dD\ap\yh_\ga\,(z\eqa y)}; 
\eay
\right\}
\eqno(\ast\ast)
$$
and hence to prove $\Yg_\ga=\Yh_\ga$ it suffices to check that 
$\yh_\ga\in\Yg_\ga$ and $\yg_\ga\in\Yh_\ga$.\vom

{\it Case 1\/}: $\ga\in\Da$. 
It follows from 2) and 3) that $\zzz_\ga\eqa \yg_\ga$
and hence $\yh_\ga\eqa y$ by 4), where $y=\vpi\ap\yg_\ga$.
Thus $\yh_\ga\in\Yg_\ga$ by $(\ast\ast)$, the first line.
On the other hand, $\zzz_\ga=\vpi^{-1}\ap \yh_\ga$ still by 4), 
so that $\yg_\ga\in\Yh_\ga$ by $(\ast\ast)$, the second line.

{\it Case 2\/}: $\ga\nin\Da$. 
Note that for a given $\ga$ 5) holds only for finitely many 
numbers $n$ by Lemma~\ref{str}\ref{str3}, so 6) holds for 
almost all $n$. 
Therefore $\yh_\ga\eqa\zzz_\ga$. 
But $\zzz_\ga=\yg_\ga$ in this case by 2). 
Thus $\yh_\ga\in\Yg_\ga$ by $(\ast\ast)$, the first line 
(with $y=\yg_\ga$).
And $\yg_\ga\in\Yh_\ga$ holds by a similar argument.
\epF{Lemma} 

\qeD{Theorem~\ref{dl1}}

\parf{The structure of the extension}
\las{F}

Here we accomplish the proof of Theorem~\ref{M}. 

\bbl
\lam{bla2}
We fix a set $G\sq \dT$, \dd\dT generic over $\bL$, during 
the course of this section.
\ebl

It will be shown that the symmetric subextension 
$\Ls[G]=\bL(\wg)$ (see Section~\ref{E}) 
satisfies Theorem~\ref{M}. 
The following is a key technical claim.

\bte 
\lam{alen}
Suppose that\/ 
$\nu<\om$, and\/ $Z\in\Ls[G]\zt Z\sq[0,\aleph_{\nu+1})$. 
Then\/ $Z\in\bL[\ens{\fg\og n}{n\le \nu}]$.
\ete 
\bpf
It follows from Lemma~\ref{redp} and 
Theorem~\ref{dl1} that there exist finite 
sets\/ $N\sq\om$ and\/ $\Ga\sq\la$ such that\/ 
$Z\in\bL[\ens{\fg\og n}{n\in N},\ens{\yg_\ga}{\ga\in\Ga}]$.
We can assume that
\ben
\aenu
\itlm{alen1}\msur
$N=\ans{0,1,2,\dots,\ka}$ for some $\ka<\om$, $\ka\ge\nu$;

\itlm{alen2}
if $\ga\ne \da$ belong to $\Ga$ and $n<\om$ satisfies 
$\ag(n,\ga)=\ag(n,\da)$ then $n\le\ka$.
\een
(Lemma~\ref{str}\ref{str3} is used to justify \ref{alen2}.) 
Define, in $\bL$, 
$$
\dT[N,\Ga]=\ens{s\in\dT}
{\bas s=N\land \abs s=\Ga\land a^s=\ag\res{(N\ti\Ga)}}.
$$

\ble
\lam{gg}
The set\/ $G[N,\Ga]=G\cap\dT[N,\Ga]$ is\/ 
\dd{\dT[N,\Ga]}generic over\/ $\bL$ and\/ $Z\in\bL[G[N,\Ga]]$.
\ele
\bpf[Lemma]
Suppose that $t\in\dT$, $N\sq\bas t$, $\Ga\sq\abs t$. 
Define the \rit{projection} $s=t[N,\Ga]\in\dT[N,\Ga]$ so that 
$q^s=q^t\res\Ga$, $a^s=a^t\res{(N\ti\Ga)}$, and $h^s$ is the 
restriction of $h^t$ to the set 
$\abs{h^t}\cap\bigcap_{n\le\ka}\aal n$.
(It is not asserted that $t\le s$.) 
Given a condition $s'\in \dT[N,\Ga]$, $s'\le s$, 
we have to accordingly find a condition $t'\le t$ such that 
$t'[N,\Ga]=s'$.

Define $t'$ as follows. 
First of all, $\bas t'=\bas t$, $\abs{t'}=\abs t$, $a^{t'}=a^t$.

Put $h^{t'}\og n=h^{s'}\og n$ for $n\in N$ 
but $h^{t'}\og n=h^{t}\og n$ for $n\in \bas t\bez N$. 

If $n\nin N$ then put $q^{t'}_\ga\og n= q^{t'}_\ga\og n$ for all 
$\ga\in\abs{t'}=\abs t$. 
If $n\in N$ and $\ga\in\abs{t'}$ then put 
$q^{t'}_\ga\og n= h^{t'}_\xi= h^{s'}_\xi$, where $\xi=a^{t'}(n,\ga)$. 
\epF{Lemma}

In continuation of the proof of the theorem, let us analyse 
$\dT[N,\Ga]$ as the forcing notion. 
It looks like the product 
$\prod_{n=0}^\ka\dH\og n\,\ti\,\dP^\Ga$: 
indeed, if $s\in\dT[N,\Ga]$ then $h^s$ can be seen as an element 
of $\prod_{n=0}^\ka\dH\og n$, $q^s$ can be seen as an element 
of $\dP^\Ga$ 
(the product of $\card\Ga$ copies of $\dP$; $\card\Ga<\om$), 
while $a^s=\ag\res{(N\ti\Ga)}$ is a constant. 
However if $n\in N$ and $\ga\in\Ga$ then 
$q^s_\ga\og n=q^t_{a^s(n,\ga)}$, hence in fact $\dT[N,\Ga]$ 
can be identified with 
$$
\textstyle
\prod_{n=0}^\ka\dH\og n\ti(\prod_{n=\ka+1}^\iy\dP\og n)^\Ga=  
\prod_{n=0}^\ka\dH\og n\ti\prod_{n=\ka+1}^\iy(\dP\og n^\Ga)\,. 
\eqno(3) 
$$
However the sets $\dP\og n$ and $\dH\og n$ as forcing notions 
are \dd{\aleph_n^+}closed, meaning that any decreasing sequence 
of length $\le\aleph_n$ has a lower bound in the same set. 
Therefore if we present $\dT[N,\Ga]$ as 
$$
\textstyle
\prod_{n=0}^\nu\dH\og n\ti   
\prod_{n=\nu+1}^\ka\dH\og n\ti\prod_{n=\ka+1}^\iy(\dP\og n^\Ga)\,,  
\eqno(4) 
$$
then it becomes clear that the second and third subproducts are 
\dd{\aleph_{n+1}^+}closed forcing notions. 
Hence, by basic results of forcing theory, 
the set $Z\sq[0,\aleph_{\nu+1})$ belongs to the 
subextension corresponding to the first subproduct 
$\prod_{n=0}^\nu\dH\og n$.
That is, $Z\in\bL[\ens{\fg\og n}{n\le \nu}]$, as required.
\epf

\bco 
\lam{Fco1}
If\/ $n<\om$ then it is true in\/ $\Ls[G]$ that\/ 
$\aleph_n$ remains a cardinal, 
the power set\/ $\cP(\aleph_n)$ is wellorderable, and\/ 
$\card(\cP(\aleph_n))=\aleph_{n+1}$.\qed
\eco

Yet cardinal preservation holds for all cardinals!

\bco 
\lam{Fle1}
Any cardinal in\/ $\bL$ remains a cardinal in\/ $\Ls[G]$. 
\eco
\bpf
Indeed we have established 
(see the proof of Theorem~\ref{alen}) 
that any set $Z\in\Ls G$, $Z\sq\bL$, belongs to a generic 
extension of $\bL$ via a forcing as in (3) in 
the proof of Theorem~\ref{alen}. 
However any such a forcing is cardinal-preserving by a simple 
cardinality argument.
\epf

To accomplish the proof of Theorem~\ref{M}, it remains to 
check that the symmetric subextension $\Ls[G]$ contains 
a surjection $\sg:2^{\omal}\onto\la$.
We define $\sg$ in $\Ls[G]$ as follows. 
If $\ga<\la$ and $z\in\Yg_\ga$ then put $\sg(z)=\ga$. 
(The definition is consistent by Lemma~\ref{cons}.)
If $z\in2^{\omal}$ does not belong to 
$\bigcup_{\ga<\la}\Yg_\ga$ then $\sg(z)=0$.
As any set $\Yg_\ga$ definitely contains $\yg_\ga$, $\sg$ 
is a surjection onto $\la$, as required.\vtm

\qeD{Theorem~\ref{M}}

\def\refname {{\Large\bf References}}

\end{document}